\documentclass[tbtags,10pt,a4paper]{amsart}

\numberwithin{equation}{section}
\allowdisplaybreaks 

\usepackage{amssymb,eucal,mathrsfs,setspace,bbm}
\usepackage[noadjust,nobreak]{cite}
\usepackage[margin=25mm, head=10mm, foot=10mm]{geometry}

\usepackage{array}
\newcolumntype{C}{>{$}c<{$}} 

\usepackage[largesc,theoremfont,nohelv]{newtxtext} 

\usepackage[libertine,cmbraces,varbb,smallerops]{newtxmath}
\begingroup
  \makeatletter
  \@for\theoremstyle:=definition,remark,plain\do{%
    \expandafter\g@addto@macro\csname th@\theoremstyle\endcsname{%
      \addtolength\thm@preskip{.5\baselineskip plus .2\baselineskip minus .2\baselineskip}
      \addtolength\thm@postskip{.5\baselineskip plus .2\baselineskip minus .2\baselineskip}
    }%
  }
\endgroup
\usepackage{enumitem}
\setitemize{leftmargin=*}       
\setenumerate{leftmargin=*,     
  label=\textnormal{\arabic*.\@}
}

\usepackage{mathtools} 
\usepackage{graphicx}

\usepackage{tikz}
\usetikzlibrary{calc}
\tikzset{%
	>=latex,
	wt/.style={circle, draw=black, fill=black, inner sep=2pt, outer sep=0pt, minimum size=5pt}, 
	nom/.style={rectangle, draw=black!20, fill=blue!20, inner sep=3pt}, 
}

\usepackage[colorlinks=true,citecolor=red,linkcolor=blue,linktoc=all]{hyperref}

\usepackage{zref-clever}                               
\zcsetup{cap,nameinlink=false}                         
\newcommand{\lcnamezcref}[1]{\zcref*[noref,nocap]{#1}} 

\theoremstyle{plain}

\theoremstyle{remark}

\AddToHook{env/proposition/begin}{\zcsetup{countertype={theorem=proposition}}}
\AddToHook{env/lemma/begin}{\zcsetup{countertype={theorem=lemma}}}
\AddToHook{env/corollary/begin}{\zcsetup{countertype={theorem=corollary}}}
\AddToHook{env/definition/begin}{\zcsetup{countertype={theorem=definition}}}
\AddToHook{env/remark/begin}{\zcsetup{countertype={theorem=remark}}}

\newcommand{\pd}{\partial}     
\newcommand{\wun}{\mathbbm{1}} 

\newcommand{\cc}{\mathsf{c}}   
\newcommand{\dd}{\mathrm{d}}   
\newcommand{\ee}{\mathrm{e}}   
\newcommand{\ii}{\mathrm{i}}   
\newcommand{\kk}{\mathsf{k}}   

\newcommand{\ahol}[1]{\overline{#1}}

\renewcommand{\cong}{\simeq}              

\renewcommand{\ge}{\geqslant} 
\renewcommand{\le}{\leqslant} 

\DeclarePairedDelimiter{\brac}{\lparen}{\rparen}   
\DeclarePairedDelimiter{\sqbrac}{\lbrack}{\rbrack} 
\DeclarePairedDelimiter{\set}{\lbrace}{\rbrace}
\newcommand{\st}{\mspace{5mu} {:} \mspace{5mu}}    
\DeclarePairedDelimiter{\abs}{\lvert}{\rvert}

\DeclarePairedDelimiterX{\comm}[2]{\lbrack}{\rbrack}{#1 , #2}  
\DeclarePairedDelimiterX{\acomm}[2]{\lbrack}{\rbrack}{#1 , #2} 
\DeclarePairedDelimiterX{\bilin}[2]{\langle}{\rangle}{#1 , #2} 
\newcommand{\no}[1]{\mathopen{:} #1 \mathclose{:}} 

\DeclarePairedDelimiter{\bra}{\langle}{\rvert}

\DeclarePairedDelimiterX{\braket}[2]{\langle}{\rangle}{#1 \delimsize\vert\mathopen{} #2} 
\DeclarePairedDelimiterX{\bracket}[3]{\langle}{\rangle}%
		{#1 \delimsize\vert\mathopen{} #2 \delimsize\vert\mathopen{} #3} 
\DeclarePairedDelimiter{\corrfn}{\langle}{\rangle} 

\DeclarePairedDelimiterXPP{\cbeta}[1]{\mathsf{B}}{\lparen}{\rparen}{}{#1}            
\DeclarePairedDelimiterXPP{\ibeta}[2]{\mathsf{B}}{\lparen}{\rparen}{}{#1;#2}         
\DeclarePairedDelimiterXPP{\twoFone}[4]{{}_2\mathsf{F}_1}{\lbrack}{\rbrack}{}{\begin{smallmatrix} #1 \ \ #2 \\ #3 \end{smallmatrix};#4} 
\DeclarePairedDelimiterXPP{\threeFtwo}[6]{{}_3\mathsf{F}_2}{\lbrack}{\rbrack}{}{\begin{smallmatrix} #1 \ \ #2 \ \ #3 \\ #4 \ \ #5 \end{smallmatrix};#6} 

\newcommand{\Ra}{\Rightarrow}

\newcommand{\lra}{\longrightarrow}

\newcommand{\ses}[3]{0 \lra #1 \lra #2 \lra #3 \lra 0}

\newcommand{\fld}[1]{\mathbbm{#1}}   
\newcommand{\alg}[1]{\mathfrak{#1}}  
\newcommand{\VOA}[1]{\mathsf{#1}}    
\newcommand{\Mod}[1]{\mathcal{#1}}   

\newcommand{\ZZ}{\fld{Z}}
\newcommand{\NN}{\fld{N}}

\newcommand{\RR}{\fld{R}}
\newcommand{\CC}{\fld{C}}

\newcommand{\SLA}[2]{\alg{#1}_{#2}}                      

\newcommand{\vac}{\Omega}                          

\newcommand{\bgsymb}{\VOA{G}}
\newcommand{\bgvoa}{\bgsymb}                       

\newcommand{\Tsing}{\widetilde{T}}                 
\newcommand{\Lsing}{\widetilde{L}}                 

\newcommand{\saff}[2]{\VOA{L}_{#1}(#2)}            

\newcommand{\specflowsymb}[1]{\sigma^{#1}}
\DeclarePairedDelimiterXPP{\specflow}[2]{\specflowsymb{#1}}{\lparen}{\rparen}{}{#2} 

\newcommand{\vacmod}{\Mod{V}}                      
\newcommand{\relmod}[1]{\Mod{W}_{[#1]}}            
\newcommand{\projmod}{\Mod{P}}                     
\newcommand{\sfvac}[1]{\specflow{#1}{\vacmod}}
\newcommand{\sfrel}[2]{\specflow{#1}{\relmod{#2}}}
\newcommand{\sfproj}[1]{\specflow{#1}{\projmod}}

\newcommand{\corrconst}[2]{\sqbrac[\big]{\begin{smallmatrix} #1 \\ #2 \end{smallmatrix}}}                    
\newcommand{\confblock}[4]{\mathscr{F}_{#1}\brac[\big]{\begin{smallmatrix} #2 \\ #3 \end{smallmatrix} ; #4}} 

\newcommand{\fuse}{\mathbin{\times}}

\newcommand{\bpz}{Belavin--Polyakov--\zam}

\newcommand{\kz}{Knizhnik--\zam}
\newcommand{\ngk}{Nahm--Gaberdiel--Kausch}

\newcommand{\wzw}{Wess--Zumino--Witten}
\newcommand{\zam}{Zamolodchikov}

\newcommand{\lhs}{left-hand side}

\newcommand{\rhs}{right-hand side}

\newcommand{\hw}{highest-weight}
\newcommand{\hwv}{\hw\ vector}
\newcommand{\hwvs}{\hwv s}
\newcommand{\hwm}{\hw\ module}

\newcommand{\chw}{conjugate highest-weight}
\newcommand{\chwv}{\chw\ vector}

\newcommand{\rhw}{relaxed highest-weight}
\newcommand{\rhwv}{\rhw\ vector}
\newcommand{\rhwvs}{\rhwv s}
\newcommand{\rhwm}{\rhw\ module}
\newcommand{\rhwms}{\rhwm s}

\newcommand{\vo}{vertex operator}
\newcommand{\voa}{\vo\ algebra}
\newcommand{\voas}{\voa s}

\newcommand{\cft}{conformal field theory}
\newcommand{\cfts}{conformal field theories}
\newcommand{\emt}{energy-momentum tensor}
\newcommand{\emts}{\emt s}
\newcommand{\ope}{operator product expansion}
\newcommand{\opes}{\ope s}
\newcommand{\qhr}{quantum hamiltonian reduction} 

\makeatletter
\renewcommand\author@andify{%
  \nxandlist {\unskip ,\penalty-1 \space\ignorespaces}%
    {\unskip {} \@@and~}%
    {\unskip \penalty-2 \space \@@and~}%
}
\makeatother

\DeclareRobustCommand{\SkipTocEntry}[5]{}

\usepackage{xpatch}
\makeatletter
\xpatchcmd{\@tocline}
          {\hfil\hbox to\@pnumwidth{\@tocpagenum{#7}}\par}
          {\ifnum#1<0\hfill\else\dotfill\fi\hbox to\@pnumwidth{\@tocpagenum{#7}}\par}
          {}{}
\makeatother

\begin{document}

\title{Bosonic ghost correlators: a case study}

\author[X~Li]{Xueting Li}
\address[Xueting Li]{
	School of Engineering \\
	Swinburne University of Technology \\
	Hawthorn, Australia, 3122.
}
\email{kittyli@swin.edu.au}

\author[D~Rajbhandari]{Damodar Rajbhandari}
\address[Damodar Rajbhandari]{
	School of Mathematics and Statistics \\
	University of Melbourne \\
	Parkville, Australia, 3010.
}
\email{drajbhandari@student.unimelb.edu.au}

\author[D~Ridout]{David Ridout}
\address[David Ridout]{
	School of Mathematics and Statistics \\
	University of Melbourne \\
	Parkville, Australia, 3010.
}
\email{david.ridout@unimelb.edu.au}

\begin{abstract}
	There has been a lot of recent work addressing the representation theory that underlies logarithmic conformal field theories.
	A full understanding of these models will however also need analytic data, in particular the correlation functions.
	Here, we explore the correlators of one of the most fundamental of all logarithmic models: the bosonic ghost system.
	In this first part, we use differential equations to show that certain correlation functions may be expressed using hypergeometric functions.
	Our main result is the consequent verification that there are four-point functions with logarithmic singularities.
	In a sequel, we will employ Coulomb gas and bootstrap methods to further refine the results presented here.
\end{abstract}

\maketitle

\markleft{X~LI, D~RAJBHANDARI AND D~RIDOUT} 

\tableofcontents

\onehalfspacing

\section{Introduction}

Conformal field theory occupies a central role in modern mathematical physics and is becoming increasingly important in areas of pure mathematics.
In particular, the chiral algebras of these theories (also known as \voas) have spurred significant advances in fields including representation theory, combinatorics, number theory, category theory and algebraic geometry.
It is therefore important to understand a range of key examples.

We can classify \cfts\ into four very broad classes based on the following two mathematical features:
\begin{itemize}
	\item A \cft\ is \emph{factorisable} if its state space is completely reducible as a module over its symmetry algebra.
	The latter is the tensor product of the chiral and antichiral \voas.
	Mathematically, this means that these \voas\ have semisimple module categories whose (Deligne) product contains the state space.
	\item A \cft\ is \emph{discrete} if its spectrum (the eigenvalues of the Virasoro zero mode) consists of isolated points.
	Mathematically, the module categories of the \voas\ are finite, meaning that there are finitely many (isomorphism classes of) simple objects.
\end{itemize}

We mention that there are models in the physics literature, see \cite{RibExa24} for example, which discuss countably infinite discrete spectra.
However, these examples have not yet received a \vo-algebraic treatment.
It would be very interesting to investigate them from this perspective.
We also note that one can further subdivide non-discrete theories by whether the fusion products of their simple objects always decompose into finitely many indecomposable objects or not.
The former situation corresponds to \emph{local finiteness} of the module category.

The four classes of \cfts\ are then as follows:
\begin{itemize}
	\item A \emph{rational} \cft\ is one that is both factorisable and discrete.
	Examples include compactified bosons, free fermions (and fermionic ghosts), Virasoro minimal models and compact \wzw\ models.
	\item A \emph{non-rational} \cft\ is one that is factorisable but not discrete.
	Examples include free bosons, which are locally finite, as well as Liouville theory and non-compact \wzw\ models, which are not.
	\item A \emph{log-rational} \cft\ is one that is discrete but not factorisable.
	Examples include symplectic fermions and the triplet models.
	\item A \emph{generic} \cft\ is one that is neither discrete nor factorisable.
	Examples include bosonic ghosts, fractional-level \wzw\ models, most supergroup \wzw\ models and most W-algebras.
\end{itemize}
Non-factorisable models are also often referred to simply as \emph{logarithmic} \cfts.
Here, we will study one of the archetypal examples of such a model: the generic \cft\ known as the bosonic ghost system (or symplectic bosons).

Originally introduced in \cite{FriCon86}, the bosonic ghosts \cft\ is a fundamental example that is essential to the construction of free-field realisations of \wzw\ models and W-algebras, see for example \cite{WakFoc86,FeiFam88,BouFre90}.
It is also a key ingredient of the definition of W-superalgebras using \qhr\ \cite{KacQua03}.
Geometrically, bosonic ghosts model the coordinates on a manifold in the chiral de~Rham complex of \cite{MalChi98} and constitute the simplest of the $\mathcal{B}_p$ algebras \cite{CreCos13} that are dual to certain $N=2$ supersymmetric 4D gauge theories under the 4D/2D correspondence of \cite{BeeInf13}.

The representation theory of the bosonic ghosts' chiral algebra was subsequently explored in \cite{LesSU202,LesLog03,RidSL208,RidFus10} using the fact that its bosonic orbifold is the $\kk=-\frac{1}{2}$ \wzw\ model associated to $\SLA{sl}{2}$.
A comprehensive analysis appeared in \cite{RidBos14}, where \rhwms\ and spectral flow were used to verify the modular invariance of characters (in the sense of the standard module formalism of \cite{CreLog13,RidVer14}).
This work also computed the fusion rules, using the standard Verlinde formula and, independently, the \ngk\ algorithm \cite{NahQua94,GabInd96}.
Recently, these latter calculations have been verified mathematically in \cite{AdaFus19,AllBos20}.

All of these works address the representation theory of the bosonic ghosts \voa.
This is an important milestone on the road to understanding a \cft.
From a physics perspective, however, it is perhaps even more important to follow up by using this information to compute correlation functions.
Despite such computations being a central feature of early studies of logarithmic \cft, see for example \cite{RozQua92,GurLog93}, it is fair to say that correlators have not received enough attention in recent times.

This paper investigates bosonic ghost correlators as a case study with which to develop tools to attack other generic \cfts.
As we shall see, the fact that the spectrum is unbounded below makes this investigation rather non-trivial.
We remark that a few ghost correlators were considered in \cite{AllBos20}, using the free-field realisation of \cite{FriCon86}.
As their focus was purely categorical, these considerations were primarily directed towards whether certain maps defined using correlators were zero or not.
Our focus is more general: we simply want to calculate some interesting correlators.
Here, we do so using standard analytic methods (solving differential equations).
In a sequel, we will explore the application of more refined methods such as the (generalisation of the) Coulomb gas approach and the conformal bootstrap.

\medskip

We now provide an outline of this paper as well as a summary of some of the key results.
\zcref{sec:bosonic-ghost-primer} first reviews the weight modules of the bosonic ghosts \voa, following \cite{RidBos14}.
This includes the classification of irreducibles in terms of spectral flows of \rhwms\ and the known fusion rules.
We also introduce the ``ghost primary'' states and fields whose correlation functions we wish to compute.

The Ward identities satisfied by $N$-point correlators are quickly reviewed in \zcref{sec:correlators}, along with their solutions for $N\le4$.
Of note here is that these identities differ slightly from the usual ones because we follow \cite{RidBos14} in choosing an asymmetric conformal structure for our ghosts: their conformal dimensions are chosen to be~$1$ and~$0$ instead of~$\frac{1}{2}$ and~$\frac{1}{2}$.
This asymmetry, originally chosen to simplify modularity computations, is convenient for our purposes because it allows us to realise the vacuum module as a submodule of a \rhwm, the latter being one of the building blocks from which the ghost theory's spectrum is built.
It is also visible in the conjugates of the primary fields, see \eqref{eq:general-2pt-function}.

An unrelated, but still important, observation is that care needs to be taken in the ghost theory with certain manipulations of correlators that are often assumed to be completely general.
In particular, one is accustomed to the vanishing of contour integrals around~$\infty$, but this need not be true here.
We analyse when this vanishing holds in \eqref{eq:big-contour} and use it to deduce our first result:~$N$-point functions of ghost primaries are~$0$ if the spectral flow indices are either all positive or all non-positive.
We therefore restrict in what follows to correlators in which a single ghost primary has a positive spectral flow index~$\ell$ while the others have index zero.

As the \emt\ of the ghost theory is a composite field, meaning that it is expressible in terms of the ghost fields, the theory admits \kz\ equations \eqref{eq:kz-method1}.
These equations are the focus of \zcref{sec:kz}.
Unfortunately, if $\ell>1$, then these equations involve correlators of descendant fields that cannot be rewritten in terms of primaries because of the issue with non-vanishing residues at~$\infty$.
This problem can however be solved by deriving a KZ equation from a carefully chosen linear combination of modes of the \emt\ and the (composite) free boson built from the ghost fields, see \eqref{eq:kz-method2}.

We remark that all of these KZ equations are recursion relations for the~$N$-point correlators, in general having no obvious initial or boundary condition.
Their solution is therefore challenging.
However, combining them with the solutions to the Ward identities leads to strong constraints on the~$N$-point functions with $N\le4$.
Just as the $2$-point functions under consideration vanish unless $\ell=1$, this combination gives our second result: the~$3$-point functions vanish unless $\ell=1$ or~$2$ (and the~$4$-point functions vanish unless $\ell=1$,~$2$ or~$3$).
Moreover, the KZ equation implies interesting recursion relations for the~$3$-point constants, see \eqref{eq:3pt-constants-l=1} and \eqref{eq:3pt-constants-l=2}.

The case $\ell=1$ is in many ways easier to analyse than $\ell>1$.
In particular, the two KZ equations are equivalent (once the Ward identities are taken into account) and involve just two correlators.
We exploit this by deriving a second ``KZ-like'' equation \eqref{eq:kz-J0}, starting not with the \emt\ of the ghosts theory but with its free boson.
For $\ell=1$, this leads to an independent equation involving the same two correlators, hence we can use it to decouple the KZ equation and thus solve it.
The solution \eqref{eq:solution-l=1} is our third result: the coordinate dependence of the~$4$-point function when $\ell=1$.
This dependence turns out to be very simple --- it is just a power law.

For $\ell>1$, we need further constraints on the correlators.
Normally, one would use null-vector relations to derive \bpz\ equations for the correlators.
However, there are no non-trivial null vectors in ghost modules (the bosonic ghosts theory is, after all, free).
Undeterred by this, we describe in \zcref{sec:bpz} our fourth result: such a BPZ equation \eqref{eq:bpz}.
Here, the loophole is that a lack of null vectors does not preclude non-trivial null relations between the modes of composite fields, in particular those of the free boson and \emt.
We moreover motivate and explain these relations as consequences of the (genuine) null-vector relations of the modules of the singlet algebra, the latter being the coset of the bosonic ghosts algebra by its free boson subalgebra \cite{RidSL210}.
Our BPZ equation corresponds to the simplest relation, which follows if we fix the charge of one of the ghost-primary fields to the value~$\frac{1}{2}$.

Restricting to correlators with one charge so-fixed, we solve our BPZ equation for the~$4$-point functions with $\ell=2$ and~$3$.
Taking into account the usual monodromy constraints, the coordinate dependence of the $\ell=3$ result is also found to be power-law in nature, see \eqref{eq:4pt-l=3}.
However, the $\ell=2$ result \eqref{eq:4pt-sol-l=2} is more interesting.
The conformal blocks in this case are hypergeometric functions and so we may tune the other charges so that the coordinate dependence of the~$4$-point function with $\ell=2$ exhibits logarithmic singularities.
This is our fifth (and in fact main) result.
Moreover, we verify that the appearance of these logarithms is in perfect agreement with the fact that, with these tuned charges, the corresponding fusion rules give modules on which the hamiltonian acts non-diagonalisably.

This achieves our main goal: to explicitly demonstrate that the bosonic ghosts \cft\ has correlation functions with logarithmic singularities.
We conclude in \zcref{sec:exploring} by using our BPZ results as initial conditions for the KZ recursion relations.
In this way, the conformal blocks found for $\ell=2,3$ and one charge fixed to $\frac{1}{2}$ are extended to charges $\frac{1}{2}, \frac{3}{2}, \frac{5}{2}, \dots$. Our sixth, and final, result is that the formulae for these blocks can be recast in a form for which this charge is completely unrestricted.
This suggests that our formulae in fact hold for all (real) charges and we test this by applying KZ-recursion to the predicted blocks with charges $-\frac{1}{2}, -\frac{3}{2}, \dots$ and show that we thereby recover the charge $\frac{1}{2}$ block.
Interestingly, for $\ell=3$, our proposed formula for general conformal blocks is hypergeometric in nature, while for $\ell=2$, we encounter generalised hypergeometric functions.

The work reported here is supported by two appendices.
The first reviews some of the basic theory of hypergeometric functions that is used.
The second sketches a derivation of an interesting identity that expresses a generalised hypergeometric function as a sum of standard ones.

\medskip

We conclude with a couple of comments concerning future work.
First, the results presented here are deliberately limited to ghost-primary correlators in which only one of the fields has a non-zero spectral flow index.
It is natural to extend these results to cover cases with more than one non-zero spectral flow and we expect that the methods used here will also apply, though the calculations are likely to be a little more involved.
This extension is relevant to a second natural extension: using bootstrap methods to fix the $3$-point constants.
Preliminary investigations indicate that this fixing will require information about more general correlators than those considered here.

A second natural question is whether our methods can be applied to bosonic ghost systems with a different central charge $\cc$.
From the point of view of the representation theory, changing the \emt\ amounts to simply ``regrading'' the modules.
One therefore expects that the correlators are also qualitatively unchanged --- there will still be logarithmic singularities in particular.
However, the quantitative features (by which we mean the conformal dimensions and $n$-point constants) will surely change.
It is not clear to us whether this data, as a function of $\cc$, can be easily inferred from that presented below.
This is an interesting problem that we leave to future work.

\addtocontents{toc}{\SkipTocEntry}
\subsection*{Acknowledgements}

We thank Jesper Jacobsen, Robert McRae, Eveliina Peltola, Sylvain Ribault and Jinwei Yang for helpful correspondence related to this work.
Some of the work reported here has previously appeared in Xueting's MSc thesis \cite{LiMSc22}.
Damodar's research is supported by an Australian Government Research Training Program (RTP) Scholarship.
David's research is supported by the Australian Research Council Discovery Project DP210101502 and an Australian Research Council Future Fellowship FT200100431.

\section{A primer on bosonic ghosts} \label{sec:bosonic-ghost-primer}

The bosonic ghosts vertex algebra was introduced in \cite{FriCon86} to describe the Faddeev--Popov ghosts that arise when gauging fermionic degrees of freedom in superstring theory.
Also known as symplectic bosons or the Weyl algebra, it is generated by two fields $\beta(z)$ and $\gamma(z)$ subject to the \opes
\begin{equation} \label{eq:opebetagamma}
	\beta(z)\beta(w) \sim 0, \quad \beta(z)\gamma(w) \sim -\frac{\wun}{z-w}, \quad \gamma(z)\gamma(w) \sim 0.
\end{equation}
This vertex algebra contains a current $J(z) = \no{\beta(z)\gamma(z)}$ and a one-parameter family of \emts.
We choose the latter to be $T(z) = -\no{\beta(z)\pd\gamma(z)}$, corresponding to central charge $\cc=2$.
The bosonic ghosts \voa\ defined by this choice will be denoted by $\bgvoa$.
Note that $\beta(z)$ and $\gamma(z)$ are now Virasoro primaries with conformal weights $1$ and $0$, respectively.
They are also charged by $J(z)$, with respective charges $+1$ and $-1$.

Given this choice of conformal structure, fields are decomposed into modes in the usual fashion:
\begin{equation} \label{eq:ope-beta-gamma}
	\beta(z) = \sum_{n\in\ZZ} \beta_n z^{-n-1}, \quad
	\gamma(z) = \sum_{n\in\ZZ} \gamma_n z^{-n}, \quad
	J(z) = \sum_{n\in\ZZ} J_n z^{-n-1}, \quad
	T(z) = \sum_{n\in\ZZ} L_n z^{-n-2}.
\end{equation}
As usual, the \opes\ imply commutation relations for the modes, such as
\begin{equation}
  \comm{\beta_m}{\gamma_n} = -\delta_{m+n=0} \wun, \quad
  \comm{J_m}{J_n} = -m \delta_{m+n=0} \wun, \quad
  \comm{L_m}{J_n} = -\frac{1}{2} m(m+1) \delta_{m+n=0} \wun - n J_{m+n}, \qquad
  m,n \in \ZZ.
\end{equation}
The unital associative algebra generated by the modes admits many automorphisms, among which are the spectral flows $\specflowsymb{\ell}$, with $\ell \in \ZZ$.
These satisfy \cite{RidBos14}
\begin{equation} \label{eq:spectral-flow}
  \specflow{\ell}{\beta_n} = \beta_{n-\ell}, \quad
  \specflow{\ell}{\gamma_n} = \gamma_{n+l}, \quad
  \specflow{\ell}{J_n} = J_n + \ell \delta_{n=0} \wun, \quad
  \specflow{\ell}{L_n} = L_n - \ell J_n - \tfrac{1}{2} \ell(\ell-1) \delta_{n=0} \wun, \qquad
  \ell,n \in \ZZ.
\end{equation}
One can also access twisted sectors of the bosonic ghosts' representation theory by allowing spectral flows with non-integer $\ell$.
We shall not consider these sectors in this work.

We turn now to the representation theory of the ghost \voa.
A weight vector in a $\bgvoa$-module is defined to be an eigenvector of $J_0$ that is also a generalised eigenvector of $L_0$.
The pair $(j,h)$ of eigenvalues of $J_0$ and $L_0$ is its weight and a weight module is one with a basis of weight vectors.
We shall refer to $j$ and $h$ as the charge and conformal weight, respectively, of the weight vector.
If the weight vector is annihilated by the $\beta_n$ and $\gamma_n$ with $n>0$, then it is a \rhwv.
One can check that \rhwvs\ are also annihilated by the $J_n$ and $L_n$ with $n>0$.
In fact, they are also annihilated by $L_0$.
Finally, a \rhwv\ that is annihilated by $\beta_0$ is a \hwv.
Such vectors are also annihilated by $J_0$.

Recall that a \hwm\ is one that is generated by a single \hwv.
Because \hwvs\ have weight $(0,0)$, $\bgvoa$ has a unique \hw\ Verma module (up to isomorphism) and it is easily verified to be irreducible.
This is of course the vacuum module $\vacmod$ and the \hwv\ is the vacuum state $\vac$.
But, there are other (non-\hw) $\bgvoa$-modules, for example those obtained from $\vacmod$ via twisting by spectral flow automorphisms.

Twisting by an automorphism such as $\specflowsymb{\ell}$ defines an invertible endofunctor on suitable module categories.
One can be concrete about this by constructing a new $\bgvoa$-module $\specflow{\ell}{\Mod{M}}$ from a known $\bgvoa$-module $\Mod{M}$.
For this, we define $\specflow{\ell}{\Mod{M}}$ to be $\Mod{M}$ as a vector space, but equipped with a different $\bgvoa$-action.
It is convenient to distinguish elements, so denote the vectors of $\specflow{\ell}{\Mod{M}}$ by $\specflow{\ell}{v}$, where $v \in \Mod{M}$.
Then, the $\bgvoa$-action is defined to be
\begin{equation} \label{eq:spectral-flow-states}
	A_n \specflow{\ell}{v} = \specflow[\big]{\ell}{\specflow{-\ell}{A_n} v}, \quad A \in \bgvoa,\ \ell,n \in \ZZ,\ v \in \Mod{M}.
\end{equation}

\begin{figure}
	\begin{tikzpicture}[xscale=0.9]
		\begin{scope}[shift={(0,0.5)}]
			\draw[->] (0,0) -- (1,0) node[right] {$j$};
			\draw[->] (0,0) -- (0,-1) node[below] {$h$};
		\end{scope}
		\begin{scope}[shift={(2,-1.75)},scale=0.4]
			\node (L) at (-3.5,1.5) {$\dots$};
			\draw[->] (L) -- node[above] {$\specflowsymb{}$} (-0.5,1.5);
			\fill[color=gray!20] (0,0) -- (2,3) -- (6,5) -- (6,0) -- cycle;
			\draw[thick] (2,3) node[label={[label distance=-3mm, text=black] above left:{$\specflow{-2}{\vac}$}}] {} -- (6,5);
			\draw[thick] (0,0) -- (2,3);
			\node[nom] at (3,-1) {$\sfvac{-2}$};
		\end{scope}
		\begin{scope}[shift={(6,-2)},scale=0.4]
			\draw[->] (-2.5,1.5) -- node[above] {$\specflowsymb{}$} (-0.5,1.5);
			\fill[color=gray!20] (0,0) -- (3,3) -- (6,3) -- (6,0) -- cycle;
			\draw[very thick, color=red] (3,3) node[label={[label distance=-3mm, text=black] above left:{$\specflow{-1}{\vac}$}}] {} -- (6,3);
			\draw[thick] (0,0) -- (3,3);
			\node[nom] at (3,-1) {$\sfvac{-1}$};
		\end{scope}
		\begin{scope}[shift={(10,-2)},scale=0.4]
			\draw[->] (-3,1) -- node[above] {$\specflowsymb{}$} (-1,1);
			\fill[color=gray!20] (0,0) -- (0,3) -- (3,3) -- (6,0) -- cycle;
			\draw[very thick, color=red] (0,3) -- (3,3) node[label={[label distance=-3mm, text=black] above right:{$\vac$}}] {};
			\draw[thick] (3,3) -- (6,0);
			\node[nom] at (3,-1) {$\vacmod$};
		\end{scope}
		\begin{scope}[shift={(16.4,-1.75)},scale=0.4]
			\draw[->] (-9,1) -- node[above] {$\specflowsymb{}$} (-7,1);
			\fill[color=gray!20] (0,0) -- (-2,3) -- (-6,5) -- (-6,0) -- cycle;
			\draw[thick] (-2,3) node[label={[label distance=-3mm, text=black] above right:{$\specflow{2}{\vac}$}}] {} -- (-6,5);
			\draw[thick] (0,0) -- (-2,3);
			\node[nom] at (-3,-1) {$\sfvac{2}$};
			\node (R) at (3.75,1.25) {$\dots$};
			\draw[->] (0.75,1.25) -- node[above] {$\specflowsymb{}$} (R);
		\end{scope}
	\end{tikzpicture}
  \caption{%
		Pictorial representation of the weights of the spectral flows of the vacuum module $\vacmod$.
		The red lines correspond to \rhwvs.%
		} \label{fig:spec-flow-hw}
\end{figure}
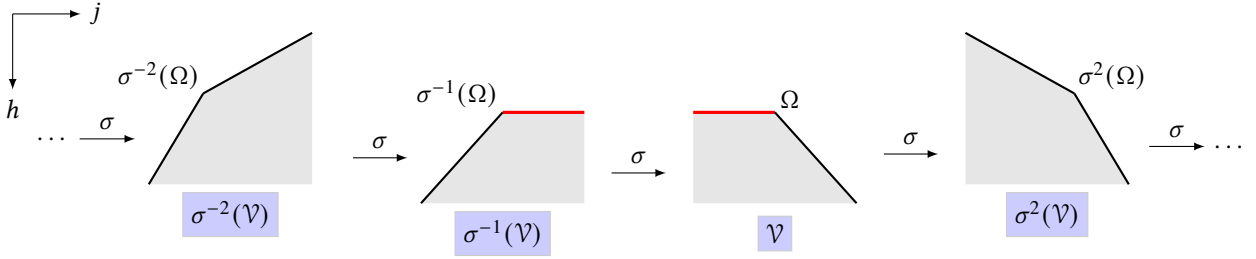

To illustrate this concrete construction, consider the action on $\sfvac{-1}$:
\begin{equation}
	\begin{aligned}
		\beta_n \specflow{-1}{\vac} &= \specflow[\big]{-1}{\specflow{}{\beta_n} \vac} = \specflow[\big]{-1}{\beta_{n-1} \vac} = 0, & \text{for all}\ n &\ge 1; \\
		\gamma_n \specflow{-1}{\vac} &= \specflow[\big]{-1}{\specflow{}{\gamma_n} \vac} = \specflow[\big]{-1}{\gamma_{n+1} \vac} = 0, & \text{for all}\ n &\ge 0.
	\end{aligned}
\end{equation}
It follows that $\specflow{-1}{\vac}$ is a \rhwv, but not a \hwv, and it is easy to check that its weight is $(1,0)$.
Being annihilated by $\gamma_0$ makes it a ``conjugate \hwv''.
Similar calculations show that $\specflow{\ell}{\vac}$ is not a \rhwv, for all $\ell \ne 0,-1$, and that the conformal weights of the states of $\sfvac{\ell}$ are not bounded below.
This is illustrated in \zcref{fig:spec-flow-hw}.

There are other \rhw\ $\bgvoa$-modules besides $\vacmod$ and $\sfvac{-1}$.
As the charge of a \rhwv\ is unconstrained, there is in fact an uncountably infinite number of them.
However, they are not parameterised by $\CC$ but by $\CC/\ZZ$ because if $v$ is a \rhwv, then (if non-zero) so are $\beta_0^n v$ and $\gamma_0^n v$ for all $n \in \NN$.
Let $\relmod{j}$ denote the \rhw\ $\bgvoa$-module for which the charges of the \rhwvs\ coincide with the coset $[j] \in \CC/\ZZ$.
These charges are thus precisely the $j+n$ with $n \in \ZZ$.
We illustrate the weights of these modules and their spectral flows in \zcref{fig:spec-flow-relax}, noting that the conformal weights of the non-trivial flows are always unbounded below.

\begin{figure}
	\begin{tikzpicture}[xscale=0.9]
		\begin{scope}[shift={(1,0)}]
			\draw[->] (0,0) -- (1,0) node[right] {$j$};
			\draw[->] (0,0) -- (0,-1) node[below] {$h$};
		\end{scope}
		\begin{scope}[shift={(4,-1.5)},scale=0.4]
			\node (L) at (-3.5,1.25) {$\dots$};
			\draw[->] (L) -- node[above] {$\specflowsymb{}$} (-0.5,1.25);
			\fill[color=gray!20] (0,0) -- (6,3) -- (6,0) -- cycle;
			\draw[thick] (0,0) -- (6,3);
			\node[nom] at (3,-1) {$\sfrel{-1}{j}$};
		\end{scope}
		\begin{scope}[shift={(8,-2)},scale=0.4]
			\draw[->] (-3,1.5) -- node[above] {$\specflowsymb{}$} (-1,1.5);
			\fill[color=gray!20] (0,0) -- (0,3) -- (6,3) -- (6,0) -- cycle;
			\draw[very thick, color=red] (0,3) -- (6,3);
			\node[nom] at (3,-1) {$\relmod{j}$};
			\draw[->] (7,1.5) -- node[above] {$\specflowsymb{}$} (9,1.5);
		\end{scope}
		\begin{scope}[shift={(14.4,-1.5)},scale=0.4]
			\fill[color=gray!20] (0,0) -- (-6,3) -- (-6,0) -- cycle;
			\draw[thick] (0,0) -- (-6,3);
			\node[nom] at (-3,-1) {$\sfrel{}{j}$};
			\node (R) at (3.75,1.25) {$\dots$};
			\draw[->] (0.75,1.25) -- node[above] {$\specflowsymb{}$} (R);
		\end{scope}
	\end{tikzpicture}
  \caption{%
		Pictorial representation of the weights of the spectral flows of a \rhwm\ $\relmod{j}$, for some $[j] \in \CC/\ZZ$.
		The red lines again correspond to \rhwvs.%
		} \label{fig:spec-flow-relax}
\end{figure}

This characterisation of the \rhwms\ by their charges turns out to be complete when the modules are irreducible, which is the case for all cosets $[j] \in \CC/\ZZ$ except $[j] = [0] = \ZZ$ \cite{RidBos14}.
In this exceptional case, there are two inequivalent \rhwms, $\relmod{0}^+$ and $\relmod{0}^-$, with the same charge coset.
They are reducible but indecomposable and may be distinguished by their submodule structures:
\begin{itemize}
	\item $\relmod{0}^+$ has a submodule isomorphic to $\vacmod$ and $\relmod{0}^+ \big/ \vacmod$ is isomorphic to $\sfvac{-1}$.
	\item $\relmod{0}^-$ has a submodule isomorphic to $\sfvac{-1}$ and $\relmod{0}^+ \big/ \sfvac{-1}$ is isomorphic to $\vacmod$.
\end{itemize}
We remark that the direct sum $\vacmod \oplus \sfvac{-1}$ also has charges in $[0] = \ZZ$.
However, it is not a \rhwm\ because it is not generated by a single \rhwv\ (it is generated by two).

For $[j] \ne [0]$, there is a unique (up to non-zero scalar multiples) \rhwv\ in $\relmod{j}$ of charge $i \in [j]$.
We shall denote it by $\phi_i$ and normalise it so that
\begin{subequations}
  \begin{equation}
    \gamma_0 \phi_i = \phi_{i-1}, \quad \text{for all}\ i \in [j].
  \end{equation}
  Since $J_0$ acts as $\gamma_0\beta_0$ on these vectors, it follows that
  \begin{equation}
    \beta_0 \phi_i = i \phi_{i+1}, \quad \text{for all}\ i \in [j].
  \end{equation}
\end{subequations}
This notation may then be extended to the \rhwvs\ of $\relmod{0}^+$, because $\phi_0 = \vac$ is a \hwv, but not to $\relmod{0}^-$.
For $\relmod{0}^-$, we shall instead denote the \rhwv\ of charge $i$ by $\widetilde{\phi}_i$, normalising it so that
\begin{equation}
  \beta_0 \widetilde{\phi}_i = \widetilde{\phi}_{i+1} \quad \text{and} \quad \gamma_0 \widetilde{\phi}_i = (i-1) \widetilde{\phi}_{i-1}
\end{equation}
($\widetilde{\phi}_1$ is then the \chwv).

Next, extend this notation to the spectral flows $\sfrel{\ell}{j}$, $\specflow{\ell}{\relmod{0}^+}$ and $\specflow{\ell}{\relmod{0}^-}$, for $\ell \in \ZZ$ and $[j] \ne [0]$, by setting $\phi_i^{\ell} = \specflow{\ell}{\phi_i}$ and $\widetilde{\phi}_i^{\ell} = \specflow{\ell}{\widetilde{\phi}_i}$ as appropriate.
It follows from \eqref{eq:spectral-flow} and \eqref{eq:spectral-flow-states} that
\begin{equation} \label{eq:def-phi}
	\begin{aligned}
		\beta_{-\ell} \phi_i^{\ell} &= i \phi_{i+1}^{\ell}, &
		\gamma_{\ell} \phi_i^{\ell} &= \phi_{i-1}^{\ell}, &
		J_0 \phi_i^{\ell} &= (i-\ell) \phi_i^{\ell}, &
		L_0 \phi_i^{\ell} &= \brac[\big]{i\ell - \tfrac{1}{2} \ell(\ell+1)} \phi_i^{\ell}, \\
		\beta_{-\ell} \widetilde{\phi}_i^{\ell} &= \widetilde{\phi}_{i+1}^{\ell}, &
		\gamma_{\ell} \widetilde{\phi}_i^{\ell} &= (i-1) \widetilde{\phi}_{i-1}^{\ell}, &
		J_0 \widetilde{\phi}_i^{\ell} &= (i-\ell) \widetilde{\phi}_i^{\ell}, &
		L_0 \widetilde{\phi}_i^{\ell} &= \brac[\big]{i\ell - \tfrac{1}{2} \ell(\ell+1)} \widetilde{\phi}_i^{\ell},
	\end{aligned}
	\qquad i \in [j],\ \ell \in \ZZ.
\end{equation}
In what follows, we shall restrict explicit attention to the $\phi_i^{\ell}$, leaving the reader to extend our results to the $\widetilde{\phi}_i^{\ell}$ as needed.
These states (and the corresponding fields) will be referred to as ghost primaries.

The fusion rules of the irreducible weight modules of the bosonic ghosts \voa\ $\bgvoa$ were first investigated (for $\cc=-1$) in \cite{RidFus10} as a corollary of a similar analysis for the fractional-level \wzw\ model $\saff{-1/2}{\SLA{sl}{2}}$.
In \cite{RidBos14}, this was revisited for $\cc=2$ using the (conjectural) standard Verlinde formula of \cite{CreLog13,RidVer14} and the \ngk\ algorithm of \cite{NahQua94,GabInd96}.
Because spectral flow indices are conserved additively under fusion \cite{LiPhy97}, it is enough to give the fusion rules of the irreducibles with no spectral flow, for example
\begin{equation} \label{eq:fusion-rule-1}
	\relmod{j_1} \fuse \relmod{j_2} \cong \relmod{j_1+j_2} \oplus \sfrel{-1}{j_1+j_2}, \qquad [j_1], [j_2], [j_1+j_2] \ne [0].
\end{equation}
This fusion rule was subsequently confirmed in \cite{AdaFus19}.

An important fact about the bosonic ghosts \voa\ $\bgvoa$ is that it admits logarithmic modules, meaning that $J_0$ acts diagonalisably but $L_0$ does not.
These modules were also first constructed, and their structures analysed, as fusion products.
In particular, when $[j_1+j_2] = [0]$ in the above fusion rule, one finds instead that
\begin{equation} \label{eq:fusion-rule-2}
	\relmod{j} \fuse \relmod{-j} \cong \sfproj{-1}, \qquad [j_1], [j_2] \ne [0],
\end{equation}
where $\projmod$ is a logarithmic $\bgvoa$-module characterised by the following non-split short exact sequences and Loewy diagram:
\begin{equation} \label{eq:proj-structure}
	\begin{gathered}
		\ses{\relmod{0}^+}{\projmod}{\specflow{}{\relmod{0}^+}}, \\
		\ses{\specflow{}{\relmod{0}^-}}{\projmod}{\relmod{0}^-},
	\end{gathered}
	\qquad\qquad
	\begin{tikzpicture}[scale=0.6,->,baseline=(l.base)]
		\node (t) at (0,2) {$\vacmod$};
		\node (l) at (-2,0) {$\sfvac{-1}$};
		\node (r) at (2,0) {$\sfvac{}$};
		\node (b) at (0,-2) {$\vacmod$};
		\draw (t) -- (l);
		\draw (t) -- (r);
		\draw (l) -- (b);
		\draw (r) -- (b);
	\end{tikzpicture}
	\ .
\end{equation}
The fusion rule \eqref{eq:fusion-rule-2} was proven in \cite{AllBos20}, where $\projmod$ was also shown to be the projective cover of the vacuum module $\vacmod$.

For completeness, we mention that \cite{AllBos20} also reported the fusion rules of $\relmod{0}^+$ and $\relmod{0}^-$:
\begin{equation} \label{eq:fusion-rule-3}
	\relmod{0}^{\pm} \fuse \relmod{0}^{\pm} \cong \relmod{0}^{\pm} \oplus \specflow{-1}{\relmod{0}^{\pm}}, \qquad
	\relmod{0}^+ \fuse \relmod{0}^- \cong \sfproj{-1}.
\end{equation}
In fact, they report the fusion rules of all (finite-length) reducible-but-indecomposable weight $\bgvoa$-modules.
We shall have no need for these in what follows.

Since $\projmod$ is a logarithmic module, the action of $L_0$ has non-trivial Jordan blocks.
Moreover, the structure \eqref{eq:proj-structure} limits the rank of these blocks to $2$ \cite{RidBos14}.
The $2$-point correlation functions of the fields corresponding to the generalised $L_0$-eigenvectors will then (generically) exhibit logarithmic singularities.
A goal of this work is to verify this feature.
To do so, we study the correlators of the ghost primary fields $\phi_i^{\ell}(z)$.
Because the fusion rule \eqref{eq:fusion-rule-2} implies that two ghost primaries are needed to generate fields in $\projmod$, we expect to see logarithmic singularities in the $4$-point functions of the $\phi_i^{\ell}(z)$.

\section{Correlation functions} \label{sec:correlators}

With the standard adjoint $L^{\dagger}_{n} = L_{-n}$ for the Virasoro modes \cite{DiFCon97}, the Ward identities give $\corrfn[\big]{\beta(z) \gamma(w)} = 0$ because the fields are Virasoro primaries of different conformal weights.
This contradicts the \ope\ \eqref{eq:ope-beta-gamma} (if we assume that the zero-point function is normalised to $1$).
The problem here is the ``asymmetric'' choice of \emt.
It is however easy to fix: the obvious adjoint on the mode algebra of $\bgvoa$ is defined by
\begin{equation} \label{eq:def-adjoint}
  \beta_n^{\dagger} = \gamma_{-n} \quad \text{and} \quad
  \gamma_n^{\dagger} = \beta_{-n}, \quad \text{hence} \quad
  J_n^{\dagger} = J_{-n} \quad \text{and} \quad
  L_n^{\dagger} = L_{-n} + n J_{-n}, \qquad n \in \ZZ.
\end{equation}
Note that $J_0$ and $L_0$ are still self-adjoint, as one should expect.

Let $\psi_i$, $i=1,\dots,N$, be a \hwv, with respect to both $J$ and $T$, of charge $j_i$ and conformal weight $h_i$:
\begin{equation} \label{eq:def-psi}
  J_n \psi_i = j_i \delta_{n=0} \psi_i, \quad
  L_n \psi_i = h_i \delta_{n=0} \psi_i, \qquad n \in \NN,\ i=1,\dots,N.
\end{equation}
Examples include the ghost primaries $\phi_i^{\ell}$ introduced in \eqref{eq:def-phi}.
The $N$-point functions of the corresponding fields then satisfy Ward identities for the modes $J_0$, $L_{-1}$, $L_0$ and $L_1$.
For the first three modes, the Ward identities have the standard form.
The identity for $L_1$ however becomes
\begin{equation} \label{eq:L1-Ward-identity}
  \sum_{i=1}^{N} \brac[\big]{w_i^2 \pd_i + 2(h_i-\tfrac{1}{2}j_i) w_i} \corrfn{\psi_1(w_1) \cdots \psi_N(w_N)} = 0.
\end{equation}

By appropriately combining the identities for $J_0$ and $L_0$, we see that the ghost Ward identities may be obtained from the standard ones by replacing each $h_i$ by $h_i-\frac{1}{2}j_i$.
This immediately implies the following modifications to the usual solutions for $N\le4$:
\begin{equation} \label{eq:Ward-identity-sols}
  \begin{aligned}
    \corrfn[\big]{\psi_1(w_1)}
      &= C_1 \delta_{h_1=0} \delta_{j=0}, \\
    \corrfn[\big]{\psi_1(w_1) \psi_2(w_2)}
      &= C_{12} \delta_{2h_1-j_1=2h_2-j_2} \delta_{j=0} w_{12}^{-2h_1+j_1}, \\
    \corrfn[\big]{\psi_1(w_1) \psi_2(w_2) \psi_3(w_3)}
      &= C_{123} \delta_{j=0} w_{12}^{-h_1-h_2+h_3-j_3} w_{13}^{-h_1+h_2-h_3-j_2} w_{23}^{h_1-h_2-h_3-j_1}, \\
    \corrfn[\big]{\psi_1(w_1) \psi_2(w_2) \psi_3(w_3) \psi_4(w_4)}
      &= \delta_{j=0} H(\eta) \:\prod_{\mathclap{1 \le a<b \le 4}}\: w_{ab}^{h-h_{ab}+j_{ab}/2}.
  \end{aligned}
\end{equation}
Here, $C_1$, $C_{12}$ and $C_{123}$ are unknown constants, $H$ is an unknown function, and we set
\begin{equation}
  w_{ab} = w_a - w_b, \quad
  \eta = \frac{w_{12} w_{34}}{w_{13} w_{24}}, \quad
  j = \sum_{i=1}^{N} j_i, \quad
  h = \frac{1}{3} \sum_{i=1}^{N} h_i, \quad
  h_{ab} = h_a + h_b \quad \text{and} \quad
  j_{ab} = j_a + j_b.
\end{equation}
The shift $h_i \mapsto h_i-\frac{1}{2}j_i$ is also required when specialising $w_i$ to $\infty$:
\begin{equation}
  \bra{\psi_i} = \lim_{w_i \to \infty} w_i^{2h_i-j_i} \bra{\vac} \psi_i(w_i).
\end{equation}
The specialised $2$-, $3$- and $4$-point functions are then
\begin{equation} \label{eq:limiting-N-point-functions}
  \begin{gathered}
    \braket{\psi_1}{\psi_2} = C_{12} \delta_{h_1-j_1=h_2} \delta_{j=0}, \quad
    \bracket{\psi_1}{\psi_2(1)}{\psi_3} = C_{123} \delta_{j=0} \\ \text{and} \quad
    \bracket{\psi_1}{\psi_2(1) \psi_3(\eta)}{\psi_4} = \delta_{j=0} \eta^{h-h_{34}+j_{34}/2} (1-\eta)^{h-h_{23}+j_{23}/2} H(\eta).
  \end{gathered}
\end{equation}

Consider now the $2$-point function with $\psi_1 = \phi_i^{\ell}$ and $\psi_2 = \phi_j^m$, where $i,j \in \CC$ and $\ell,m \in \ZZ$.
From \eqref{eq:def-phi} and \eqref{eq:Ward-identity-sols}, this vanishes unless $i+j=\ell+m$ and $(2\ell-1)i-\ell^2 = (2m-1)j-m^2$.
Solving these equations gives two solutions: either $i+j=1$ and $\ell+m=1$ or $i=j$ and $\ell+m=2i$.
The latter solution must be discarded as $2i$ need not be an integer, leaving us with the $2$-point function
\begin{equation} \label{eq:general-2pt-function}
	\corrfn[\big]{\phi_i^{\ell}(w_1) \phi_j^m(w_2)} = \frac{\corrconst{\ell&m}{i&j} \delta_{i+j=1}^{\ell+m=1}}{w_{12}^{(2\ell-1)i-\ell^2}},
\end{equation}
where $\corrconst{\ell&m}{i&j}$ denotes the $2$-point constant.
Consequently, we may choose the conjugate of the ghost primary field $\phi_i^{\ell}(z)$ to be $\phi_{1-i}^{1-\ell}(z)$, hence the conjugate of the $\bgvoa$-module $\sfrel{\ell}{i}$ is $\sfrel{1-\ell}{-i}$.
This agrees with the notion of conjugation introduced in \cite{RidBos14} as a twist by an automorphism.
It also agrees with the fusion rule \eqref{eq:fusion-rule-2} because it implies that
\begin{equation}
	\sfrel{\ell}{i} \fuse \sfrel{1-\ell}{-i} \cong \specflow{}{\relmod{i} \fuse \relmod{-i}} \cong \projmod
\end{equation}
and $\projmod$ is the projective cover (and injective hull) of the vacuum module $\vacmod$.

We can also deduce constraints on general $N$-point functions.
For this, we first recall a standard trick \cite{DiFCon97} of \cft.
Let $A(z)$ be a field of conformal weight $h_A$ and consider a correlation function involving a descendant field of the form $(A_n \psi)(w)$ (where $\psi(w)$ is an arbitrary field).
We start from
\begin{equation}
  \corrfn[\big]{\cdots (A_n \psi)(w) \cdots}
  = \oint_w \corrfn[\big]{\cdots A(z) \psi(w) \cdots} (z-w)^{n+h_A-1} \frac{\dd z}{2\pi\ii}
\end{equation}
and then write the contour around $w$ as a big contour around all insertion points (so effectively around $\infty$) minus contours around each insertion point except $w$.
However, this trick is only effective if the big contour contributes nothing.
But,
\begin{equation} \label{eq:big-contour}
  \begin{split}
    \oint_{\infty} \bra{\vac} A(z) (z-w)^{n+h_A-1} \frac{\dd z}{2\pi\ii}
    &= -\oint_0 \bra{\vac} A(y^{-1}) (y^{-1}-w)^{n+h_A-1} \frac{\dd y}{-2\pi\ii y^2} \\
    &= \sum_{r\in\ZZ} \bra{A^{\dag}_{-r} \vac} \oint_0 y^{r-n-1} (1-wy)^{n+h_A-1} \frac{\dd y}{2\pi\ii}.
  \end{split}
\end{equation}
The bra is zero unless $r$ is sufficiently large, say $r \ge h_{A^{\dag}}$.
Moreover, if in addition $n < h_{A^{\dag}}$, then $r-n-1 > -1$ and the residue at $0$ vanishes.
We conclude that the big contour contributes nothing whenever $n < h_{A^{\dag}}$.

Consider now $A = \gamma$, so that the big contour \eqref{eq:big-contour} vanishes for all $n < h_{\beta} = 1$.
Then, the relaxed primary $\phi_i(w) = (\gamma_0 \phi_{i+1})(w)$ is a descendant to which the above trick applies.
If another field in the correlator is also a relaxed primary, $\phi_j(x)$ say, then its contribution after applying the trick is
\begin{equation}
  -\oint_x \corrfn[\big]{\cdots \gamma(z) \phi_j(x) \cdots} (z-w)^{-1} \frac{\dd z}{2\pi\ii}
  = 0,
\end{equation}
since $\gamma(z) \phi_j(x) = \no{\gamma(z) \phi_j(x)}$ is regular at $z=x$.
It follows that every correlator of relaxed primaries vanishes (provided that at least one is a $\gamma_0$-descendant).

This generalises: for any $\ell\le0$, consider the ghost primary field $\phi_i^{\ell}(w) = (\gamma_{\ell} \phi_{i+1}^{\ell})(w)$ as a descendant of the ghost primary field $\phi_{i+1}^{\ell}(w)$.
If $\phi_j^m(x)$, with $m\le0$, is also in the correlator, then its contribution vanishes because its \ope\ with $\gamma(z)$ is regular.
We conclude that every correlator of spectral flows of ghost primaries vanishes when the spectral flow indices are all non-positive.
A similar argument with $\beta_{-\ell}$, for $\ell\ge1$, shows that correlators also vanish when the spectral flow indices are all positive.

\section{\kz\ equations} \label{sec:kz}

To get more detailed information about the correlation functions, we need to derive more constraints.
As $T(z)$ is not a generating field of $\bgvoa$, we can use the explicit form of $L_{-1}$ to derive \kz\ (KZ) equations.
The precise form of this equation will depend on the spectral flow indices of the fields appearing in our $N$-point function.
Given the results of the previous \lcnamezcref{sec:correlators}, a simple (non-trivial) choice is to give the first $N-1$ fields index $0$ and the last index $\ell>0$:
\begin{equation} \label{eq:standard-correlator}
	\corrfn[\big]{\phi_{j_1}(w_1) \cdots \phi_{j_{N-1}}(w_{N-1}) \phi_{j_N}^{\ell}(w_N)}.
\end{equation}

Acting with $L_{-1}$ on the $i$-th field, where $i \ne N$, now gives $L_{-1} \phi_{j_i} = -\gamma_{-1} \beta_0 \phi_{j_i} = -j_i \gamma_{-1} \phi_{j_i+1}$.
Manipulating in the usual way, we arrive at the KZ equation
\begin{multline} \label{eq:kz-method1}
  \pd_i \corrfn[\big]{\phi_{j_1}(w_1) \cdots \phi_{j_{N-1}}(w_{N-1}) \phi_{j_N}^{\ell}(w_N)} \\
  = j_i \sum_{k=1}^{\ell} \frac{k}{w_{iN}^{k+1}} \corrfn[\big]{\phi_{j_1}(w_1) \cdots \phi_{j_i+1}(w_i) \cdots \phi_{j_{N-1}}(w_{N-1}) (\gamma_k \phi_{j_N}^{\ell})(w_N)},
  \quad i \ne N.
\end{multline}
Unfortunately, if $\ell>1$, then there are correlators on the \rhs\ that involve fields $(\gamma_k \phi_{j_N}^{\ell})(w_N)$ that are not ghost primaries.
Even worse, because $h_{\gamma^{\dag}} = h_{\beta} = 1$ and $k\ge1$, we cannot simplify these correlators with the usual tricks without attracting a contribution from a residue at $\infty$.

Alternatively, we may act with $L_{-1}$ on the $N$-th field, resulting in
\begin{equation}
  \begin{split}
    \pd_N \corrfn[\big]{\phi_{j_1}(w_1) \cdots \phi_{j_{N-1}}(w_{N-1}) \phi_{j_N}^{\ell}(w_N)}
    &= \ell \corrfn[\big]{\phi_{j_1}(w_1) \cdots \phi_{j_{N-1}}(w_{N-1}) (\beta_{-\ell-1} \phi_{j_N-1}^{\ell})(w_N)} \\
    &\quad + j_N (\ell-1) \corrfn[\big]{\phi_{j_1}(w_1) \cdots \phi_{j_{N-1}}(w_{N-1}) (\gamma_{\ell-1}\phi_{j_N+1}^{\ell})(w_N)}.
  \end{split}
\end{equation}
Since $\ell>0$ and $h_{\beta^{\dag}} = h_{\gamma} = 0$, the first term may be simplified in the usual way without interference from residues at $\infty$.
However, the second cannot be simplified unless $\ell=1$.
Happily, we may instead cancel this second term by noting that $\brac[\big]{L_{-1} - (\ell-1) J_{-1}} \phi_{j_N}^{\ell} = \beta_{-\ell-1} \phi_{j_N-1}^{\ell}$, giving
\begin{equation} \label{eq:kz-method2}
  \begin{split}
    \brac[\bigg]{\pd_N + (\ell-1) \sum_{i=1}^{N-1} \frac{j_i}{w_{iN}}} &\corrfn[\big]{\phi_{j_1}(w_1) \cdots \phi_{j_{N-1}}(w_{N-1}) \phi_{j_N}^{\ell}(w_N)} \\
    &= -\sum_{i=1}^{N-1} \frac{j_i}{w_{iN}^{\ell+1}} \corrfn[\big]{\phi_{j_1}(w_1) \cdots \phi_{j_i+1}(w_i) \cdots \phi_{j_{N-1}}(w_{N-1}) \phi_{j_N-1}^{\ell}(w_N)}.
  \end{split}
\end{equation}
Unlike \eqref{eq:kz-method1}, this KZ equation only involves spectral flows of relaxed primary fields.

When $N=1$, the KZ equation \eqref{eq:kz-method2} reduces to the $L_{-1}$ Ward identity.
We get a little more information when $N=2$ by substituting instead the following version of \eqref{eq:general-2pt-function}:
\begin{equation} \label{eq:special-2pt-function}
	\corrfn[\big]{\phi_{j_1}(w_1) \phi_{j_2}^1(w_2)} = \corrconst{0&1}{j_1&j_2} \delta_{j_1+j_2=1} w_{12}^{j_1}.
\end{equation}
The result is that this is a solution provided that the $2$-point constants satisfy
\begin{equation}
	\corrconst{0&1}{j_1-1&j_2+1} = \corrconst{0&1}{j_1&j_2}, \quad j_1+j_2=1.
\end{equation}

The normalisation of the ghost primary fields $\phi_i^{\ell}(z)$, $i \in [j]$, is fixed by \eqref{eq:def-phi}, once we choose it for an arbitrary choice of $i \in [j]$ (with the proviso that our choice satisfies $i>0$ if $[j]=[0]$).
For each $[j] \in \CC/\ZZ$, we may therefore normalise the chosen $\phi_i^{\ell}$, $i \in [j]$, so that $\corrconst{0&1}{i&1-i} = 1$.
Consequently, the $2$-point constants are actually constant (independent of $j$):
\begin{equation} \label{eq:2pt-constants}
	\corrconst{0&1}{j&1-j} = 1, \quad j \in \CC.
\end{equation}
(Here, we reject the possibility that any of the $\corrconst{0&1}{j&1-j}$ vanishes.
This would mean that the corresponding field has no conjugate, contradicting the physical consistency requirement that the $2$-point functions are non-degenerate.)

The calculations are somewhat more cumbersome when $N=3$, but also more rewarding.
Substituting
\begin{equation} \label{eq:special-3pt-function}
	\corrfn[\big]{\phi_{j_1}(w_1) \phi_{j_2}(w_2) \phi_{j_3}^{\ell}(w_3)} = \corrconst{0&0&\ell}{j_1&j_2&j_3} \delta_{j_1+j_2+j_3=\ell} w_{12}^{(j_3-\ell/2)(\ell-1)} w_{13}^{-j_2 - (j_3-(\ell+1)/2)\ell} w_{23}^{-j_1 - (j_3-(\ell+1)/2)\ell}
\end{equation}
into \eqref{eq:kz-method2} gives the following polynomial equation in $w_{13}$ and $w_{23}$ (valid for $j_1+j_2+j_3=\ell$):
\begin{multline} \label{eq:3pt-function-constraint}
	\corrconst{0&0&\ell}{j_1&j_2&j_3} w_{12}^{\ell-1} \sqbrac[\Big]{\brac[\big]{j_1 + j_2(\ell-1) + (j_3-\tfrac{\ell+1}{2})\ell} w_{13} + \brac[\big]{j_1(\ell-1) + j_2 + (j_3-\tfrac{\ell+1}{2})\ell} w_{23}} \\
	= -j_2 \corrconst{0&0&\ell}{j_1&j_2+1&j_3-1} w_{13}^{\ell} - j_1 \corrconst{0&0&\ell}{j_1+1&j_2&j_3-1} w_{23}^{\ell}.
\end{multline}

For $\ell=1$, this reduces to the following identity of $3$-point constants:
\begin{equation} \label{eq:3pt-constants-l=1}
	\corrconst{0&0&1}{j_1+1&j_2&j_3-1} = \corrconst{0&0&1}{j_1&j_2&j_3} = \corrconst{0&0&1}{j_1&j_2+1&j_3-1}.
\end{equation}
(Of course, we cannot infer from this that $3$-point constants are also charge-independent when $\ell=1$.)
For $\ell=2$ however, we instead find that the $3$-point constants satisfy
\begin{equation} \label{eq:3pt-constants-l=2}
	j_1 \corrconst{0&0&2}{j_1+1&j_2&j_3-1} = (j_3-1) \corrconst{0&0&2}{j_1&j_2&j_3} \quad \text{and} \quad
	j_2 \corrconst{0&0&2}{j_1&j_2+1&j_3-1} = -(j_3-1) \corrconst{0&0&2}{j_1&j_2&j_3}.
\end{equation}
Interestingly, when $j_1$, $j_2$ or $j_3$ is an integer, hence the corresponding modules are reducible, some of these $3$-point constants vanish.
Finally, the $3$-point constants with $\ell>2$ always vanish because the factor $w_{12}^{\ell-1} = (w_{13}-w_{23})^{\ell-1}$ of the \lhs\ of \eqref{eq:3pt-function-constraint} yields cross-terms that cannot be accommodated by the \rhs.

The implication for the (generic) fusion rules is easily unravelled.
Since the $3$-point function \eqref{eq:special-3pt-function} vanishes unless $\ell=1$ or $2$, the \ope\ of $\phi_{j_1}(w_1)$ and $\phi_{j_2}(w_2)$ must produce the field conjugate to $\phi_{j_3}^{\ell}(w_3)$, which is $\phi_{1-j_3}^{1-\ell}(w_2)$ according to \eqref{eq:general-2pt-function}.
Since $j_1+j_2+j_3 = \ell$, this conjugate field is $\phi_{1-j_3}(w_2) = \phi_{j_1+j_2}(w_2)$, when $\ell=1$, and $\phi_{1-j_3}^{-1}(w_2) = \phi_{j_1+j_2-1}^{-1}(w_2)$, when $\ell=2$.
We therefore deduce that
\begin{equation} \label{eq:fusion-rule-1-from-KZ}
	\relmod{j_1} \fuse \relmod{j_2} \subseteq \relmod{j_1+j_2} \oplus \specflow{-1}{\relmod{j_1+j_2}}, \quad [j_1], [j_2], [j_1+j_2] \ne [0],
\end{equation}
with equality if $\corrconst{0&0&1}{j_1&j_2&1-j_1-j_2}$ and $\corrconst{0&0&2}{j_1&j_2&2-j_1-j_2}$ are both non-zero.
This is of course consistent with \eqref{eq:fusion-rule-1}.

We next turn to the $4$-point function
\begin{equation} \label{eq:special-4pt-function}
	\corrfn[\big]{\phi_{j_1}(w_1) \phi_{j_2}(w_2) \phi_{j_3}(w_3) \phi_{j_4}^{\ell}(w_4)}
	= \delta_{j_1+j_2+j_3+j_4=\ell} H\corrconst{0&0&0&\ell}{j_1&j_2&j_3&j_4}(\eta)
	\:\prod_{\mathclap{1 \le a<b \le 3}}\: w_{ab}^{h_4^{\ell}/3+j_{ab}/2} \:\prod_{\mathclap{1 \le c \le 3}}\: w_{c4}^{-2h_4^{\ell}/3+(j_{c4}-\ell)/2},
\end{equation}
where $H\corrconst{0&0&0&\ell}{j_1&j_2&j_3&j_4}$ is an unknown function of $\eta = \frac{w_{12}w_{34}}{w_{13}w_{24}}$.
While the KZ equation \eqref{eq:kz-method2} for $N=2$ and $3$ resulted in infinite-order recursion relations among the $2$- and $3$-point constants, respectively, substituting \eqref{eq:special-4pt-function} into \eqref{eq:kz-method2} will result in an infinite system of coupled differential equations.

Rather than try to solve this system immediately, let us first consider the consequences of the fusion constraints \eqref{eq:fusion-rule-1-from-KZ} for the $4$-point functions \eqref{eq:special-4pt-function}.
As noted above, the only ghost primaries that may appear in the \ope\ of $\phi_{j_1}(w_1)$ and $\phi_{j_2}(w_2)$ are $\phi_{j_1+j_2}(w_2)$ and $\phi_{j_1+j_2-1}^{-1}(w_2)$.
Similarly, the only ghost primaries that may appear when expanding $\phi_{j_3}(w_3) \phi_{j_4}^{\ell}(w_4)$ are $\phi_{j_3+j_4}^{\ell}(w_4) = \phi_{\ell-j_1-j_2}^{\ell}(w_4)$ and $\phi_{j_3+j_4-1}^{\ell-1}(w_4) = \phi_{\ell-1-j_1-j_2}^{\ell-1}(w_4)$.
(This follows from the general identity $\Mod{M} \fuse \specflow{\ell}{\Mod{N}} \cong \specflow{\ell}{\Mod{M} \fuse \Mod{N}}$ connecting fusion and spectral flow, see \cite{LiPhy97}).
The $4$-point function must vanish unless these two \opes\ produce conjugate fields.
But the conjugate fields of $\phi_{j_1+j_2}(w_2)$ and $\phi_{j_1+j_2-1}^{-1}(w_2)$ are $\phi_{1-j_1-j_2}^1(w_2)$ and $\phi_{2-j_1-j_2}^2(w_2)$, respectively.
We therefore conclude that the $4$-point function \eqref{eq:special-4pt-function} vanishes unless $\ell=1$, $2$ or $3$.

Consider the simplest non-vanishing possibility: $\ell=1$.
In this case, the KZ equations \eqref{eq:kz-method1} and \eqref{eq:kz-method2} involve only ghost primary fields.
However, they are not independent --- summing the first from $i=1$ to $N-1$, then adding the second, recovers the $L_{-1}$ Ward identity.
We therefore focus on \eqref{eq:kz-method1} with $i=3$, because it may be simplified by sending $w_1$ to $\infty$, $w_2$ to $1$ and $w_4$ to $0$.
Then, $w_3$ is sent to $\eta$ and $\pd_3$ to $\pd_{\eta}$, resulting in the simplified KZ equation
\begin{equation} \label{eq:kz-l=1}
	\pd_{\eta} \bracket{\phi_{j_1}}{\phi_{j_2}(1) \phi_{j_3}(\eta)}{\phi_{j_4}^1}
	= \frac{j_3}{\eta^2} \bracket{\phi_{j_1}}{\phi_{j_2}(1) \phi_{j_3+1}(\eta)}{\phi_{j_4-1}^1}.
\end{equation}
This infinite system of ordinary differential equations connects two correlators that differ only in the charges of the last two fields.
Unfortunately, we have no boundary conditions for these charges.

Because there are only two correlators in \eqref{eq:kz-l=1}, this system may be decoupled by deriving a second independent system involving the same correlators.
Such a second system is indeed available because we can deduce a second family of KZ equations from the fact that $J(z)$ is not a generating field of $\bgvoa$.
Indeed, inserting $j_i \phi_{j_i} = J_0 \phi_{j_i} = \gamma_0 \beta_0 \phi_{j_i}$, $i \ne N$, into the correlator \eqref{eq:standard-correlator}, we deduce in the usual way the following algebraic KZ equation:
\begin{multline} \label{eq:kz-J0}
	\corrfn[\big]{\phi_{j_1}(w_1) \cdots \phi_{j_{N-1}}(w_{N-1}) \phi_{j_N}^{\ell}(w_N)} \\
	= \sum_{k=1}^{\ell} \frac{\corrfn[\big]{\phi_{j_1}(w_1) \cdots \phi_{j_i+1}(w_i) \cdots \phi_{j_{N-1}}(w_{N-1}) (\gamma_k\phi_{j_N}^{\ell})(w_N)}}{w_{iN}^k}, \quad i \ne N.
\end{multline}

Setting $i=3$, $\ell=1$ and $N=4$ in \eqref{eq:kz-J0}, then sending $w_1$ to $\infty$, $w_2$ to $1$ and $w_4$ to $0$, we obtain
\begin{equation} \label{eq:kz-J0-l=1}
	\bracket{\phi_{j_1}}{\phi_{j_2}(1) \phi_{j_3}(\eta)}{\phi_{j_4}^1}
	= \frac{1}{\eta} \bracket{\phi_{j_1}}{\phi_{j_2}(1) \phi_{j_3+1}(\eta)}{\phi_{j_4-1}^1}.
\end{equation}
Substituting into \eqref{eq:kz-l=1} then decouples the system of differential equations as promised:
\begin{equation} \label{eq:kz-l=1-decoupled}
	\pd_{\eta} \bracket{\phi_{j_1}}{\phi_{j_2}(1) \phi_{j_3}(\eta)}{\phi_{j_4}^1}
	= \frac{j_3}{\eta} \bracket{\phi_{j_1}}{\phi_{j_2}(1) \phi_{j_3}(\eta)}{\phi_{j_4}^1}.
\end{equation}
The solutions are now trivial to find:
\begin{equation} \label{eq:solution-l=1}
	\bracket{\phi_{j_1}}{\phi_{j_2}(1) \phi_{j_3}(\eta)}{\phi_{j_4}^1} = \delta_{j_1+j_2+j_3+j_4=1} \corrconst{0&0&0&1}{j_1&j_2&j_3&j_4} \eta^{j_3}.
\end{equation}
Unfortunately, they do not exhibit the logarithmic singularities that we seek.

We can nevertheless extract a little more information from these solutions, specifically about the $4$-point constants.
For this, note that the counterparts of \eqref{eq:kz-J0-l=1} for $i=1,2$ are even simpler:
\begin{equation} \label{eq:kz-J0-l=1-i=1,2}
	\bracket{\phi_{j_1}}{\phi_{j_2}(1) \phi_{j_3}(\eta)}{\phi_{j_4}^1}
	= \bracket{\phi_{j_1+1}}{\phi_{j_2}(1) \phi_{j_3}(\eta)}{\phi_{j_4-1}^1}
	= \bracket{\phi_{j_1}}{\phi_{j_2+1}(1) \phi_{j_3}(\eta)}{\phi_{j_4-1}^1}.
\end{equation}
Substituting the solution \eqref{eq:solution-l=1} back into these KZ equations, as well as into \eqref{eq:kz-J0-l=1} itself, we deduce (by-now-expected) infinite-order recursion relations for the $4$-point constants:
\begin{equation} \label{eq:4pt-constants-l=1}
	\corrconst{0&0&0&1}{j_1&j_2&j_3&j_4} = \corrconst{0&0&0&1}{j_1+1&j_2&j_3&j_4-1}
	= \corrconst{0&0&0&1}{j_1&j_2+1&j_3&j_4-1} = \corrconst{0&0&0&1}{j_1&j_2&j_3+1&j_4-1}.
\end{equation}
(Again, this does not prove that the $4$-point constants for $\ell=1$ are charge-independent.)

Unfortunately, the KZ equations derived above do not decouple for $\ell>1$.
We can see this by observing that for $i=3$ say, \eqref{eq:kz-method1} and \eqref{eq:kz-J0} reduce to two equations relating $\ell+1$ independent correlators.
Moreover, all but two of these correlators involve fields that are not ghost primaries.
When $\ell=2$, one can eliminate the single correlator with a non-primary field, but the result is a single differential equation involving two primary correlators.
The situation when $\ell=3$ is even worse.
We need a different approach.

\section{A \bpz\ equation} \label{sec:bpz}

A \bpz\ (BPZ) equation is a constraint on a correlation function that is deduced from the vanishing of a null vector.
More precisely, such equations arise when the physical module, for instance an irreducible one, is a proper quotient of another module, for instance a Verma module.
However, none of the irreducible $\bgvoa$-modules we are considering have (non-zero) null vectors --- they are all (spectral flows of relaxed) Verma modules and not proper quotients.
It follows that there are no (non-zero) null states from which to deduce BPZ equations.

We will nevertheless derive such an equation.
To achieve this, we exploit a loophole in the above argument.
Observe that if we try to express a putative null state in a $\bgvoa$-module using $\beta$- and $\gamma$-modes, then the terms will necessarily cancel one another and the null state will vanish identically.
However, this cancellation need not occur if we look for null states involving only $J$- and $L$-modes.
(Of course, if we rewrite the $J$- and $L$-modes in terms of $\beta$- and $\gamma$-modes, then cancellation is guaranteed.
The perhaps surprising point is that we can deduce a non-trivial BPZ equation without such a rewriting.)

To motivate this idea, recall that the coset of the bosonic ghost system by the free boson (generated by) $J$ is the singlet algebra of central charge $\cc=-2$ \cite{RidSL210}.
The latter algebra's irreducible modules are all direct sums of irreducible Virasoro modules \cite{AdaCla02}, so null vectors will exist whenever the conformal weight belongs to the extended Kac table for $\cc=-2$.
This is the set
\begin{equation} \label{eq:ext-Kac-table}
	\set[\big]{\tfrac{1}{8}s(s-2) \st s \in \NN} = \set[\big]{-\tfrac{1}{8}, 0, \tfrac{3}{8}, 1, \tfrac{15}{8}, 3, \tfrac{35}{8}, 6, \dots}.
\end{equation}

The singlet \emt\ $\Tsing(z)$ may be identified in $\bgvoa$ by simply imposing the condition that it commutes with $J(w)$.
It is easy to verify that the unique solution is
\begin{equation}
	\Tsing(z) = T(z) + \frac{1}{2} \no{J(z)J(z)} + \frac{1}{2} \pd J(z).
\end{equation}
Indeed, a part of this verification shows that $\Tsing(z)$ has central charge $-2$.
The action of the singlet zero mode $\Lsing_0$ on a ghost primary $\phi_j$ is then
\begin{equation}
	\Lsing_0 \phi_j = \frac{1}{2}j(j-1) \phi_j.
\end{equation}
Comparing with \eqref{eq:ext-Kac-table}, there is a match if $j=\frac{1}{2}s$ or $1-\frac{1}{2}s$.
We conclude that there is a null vector descended from $\phi_j$ using $\Lsing$-modes if and only if $j \in \frac{1}{2}\ZZ$.

Let us see this in action.
For $j=0$, the null vector has grade $1$, hence it is
\begin{equation}
	\Lsing_{-1} \phi_0 = L_{-1} \phi_0 + J_{-1}J_0 \phi_0 + \frac{1}{2} (\pd J)_{-1} \phi_0 = L_{-1} \phi_0.
\end{equation}
Recalling that $\beta_0 \phi_0 = 0$, so $\phi_0$ may be identified with the vacuum $\vac$ (under $\vacmod \subset \relmod{0}$), this null vector just confirms the translation invariance of the vacuum.
It does not give us a new BPZ equation.
On the other hand, we can verify that it vanishes identically when expressed in terms of $\beta$- and $\gamma$-modes: $L_{-1} \phi_0 = -\beta_{-1}(\pd\gamma)_0 \phi_0 - (\pd\gamma)_{-1}\beta_0 \phi_0 = 0$.

The first useful null vector has $j=\frac{1}{2}$ and grade $2$:
\begin{equation} \label{eq:sing-vect}
	\chi = (\Lsing_{-1}^2 - \tfrac{1}{2} \Lsing_{-2}) \phi_{1/2} = (L_{-1}^2 - \tfrac{1}{2} L_{-2} + J_{-1}L_{-1}) \phi_{1/2}.
\end{equation}
It likewise can be checked to vanish identically when expressed in terms of $\beta$- and $\gamma$-modes, an exercise we leave to the reader.
Nevertheless, inserting it into a correlator and applying the usual methods leads to the following BPZ equations:
\begin{equation} \label{eq:bpz}
	\sqbrac[\bigg]{\pd_i^2 + \sum_{k \ne i} \frac{j_k}{w_{ik}}\pd_i - \frac{1}{2} \sum_{k \ne i} \brac[\bigg]{\frac{1}{w_{ik}}\pd_k + \frac{h_k}{w_{ik}^2}}} \corrfn[\big]{\psi_1(w_1) \cdots \psi_N(w_N)} = 0, \quad i=1,\dots,N.
\end{equation}
Here, the $\psi_1, \dots, \psi_N$ are \hwvs\ with respect to $J$ and $T$, as in \eqref{eq:def-psi}, and $\psi_i$ is fixed to be $\phi_{1/2}$.

Restricting now to correlators of the form \eqref{eq:standard-correlator}, it is easy to check that the $2$-point function \eqref{eq:special-2pt-function} is a solution of our BPZ equation \eqref{eq:bpz} (with $i=1$, $j_1=\frac{1}{2}$ and $N=2$):
\begin{equation}
	\sqbrac[\bigg]{\pd_1^2 + \frac{j_2}{w_{12}}\pd_1 - \frac{1}{2w_{12}}\pd_2 - \frac{j_2-1}{2w_{12}^2}} \corrfn[\big]{\phi_{1/2}(w_1) \phi_{j_2}^1(w_2)} = 0.
\end{equation}
The same is true for the $3$-point function \eqref{eq:special-3pt-function}, but only when $\ell=1$ or $2$:
\begin{multline} \label{eq:3pt-BPZ}
	\sqbrac[\bigg]{\pd_1^2 + \brac[\Big]{\frac{j_2}{w_{12}} + \frac{j_3-\ell}{w_{13}}}\pd_1 - \frac{1}{2w_{12}}\pd_2 - \frac{1}{2w_{13}}\pd_3 - \frac{j_3\ell - \frac{1}{2}\ell(\ell+1)}{2w_{13}^2}} \corrfn[\big]{\phi_{1/2}(w_1) \phi_{j_2}(w_2) \phi_{j_3}^{\ell}(w_3)} \\
	= (\ell-1)(\ell-2) \brac[\big]{j_2-\tfrac{1}{2}\ell} \brac[\big]{j_2-\tfrac{1}{2}(\ell-1)} w_{12}^{-j_2(\ell-1)+(\ell^2-2\ell-3)/2} w_{13}^{j_2(\ell-1)-(\ell^2-2\ell+4)/2} w_{23}^{j_2\ell-(\ell^2-2\ell-3)/2}.
\end{multline}
For $j_1=\frac{1}{2}$, the BPZ equation thus reproduces the vanishing of the $3$-point functions with $\ell\ne1,2$, and hence the fusion rule \eqref{eq:fusion-rule-1-from-KZ}, without invoking the KZ equation.

Consider next the $4$-point function \eqref{eq:special-4pt-function} with $j_3=\frac{1}{2}$ and $\ell=1$.
Substituting it into the ($i=3$, $N=4$) BPZ equation \eqref{eq:bpz} results, after some work, in a second-order ordinary differential equation for the unknown function $H(\eta) = H\corrconst{0&0&0&1}{j_1&j_2&j_3&j_4}(\eta)$.
The precise form of this equation is rather long, so we omit it and instead report the result of substituting the solution back into the specialised correlator with $w_1 \to \infty$, $w_2 \to 1$ and $w_4 \to 0$:
\begin{equation} \label{eq:4pt-sol-l=1}
	\bracket[\big]{\phi_{j_1}}{\phi_{j_2}(1) \phi_{1/2}(\eta)}{\phi_{j_4}^1}
	= \delta_{j_1+j_2+j_4=1/2} \eta^{1/2} \brac[\Big]{\alpha_1 + \alpha_2\ibeta[\big]{-j_4+\tfrac{1}{2}, -j_2+\tfrac{1}{2}}{\eta}}.
\end{equation}
Here, $\alpha_1$ and $\alpha_2$ are constants (depending on the charges, but independent of $\eta$) and $\ibeta{a,b}{z}$ is the incomplete beta function \eqref{eq:def-Beta}.
The conformal blocks in this case are thus $\eta^{1/2}$ and $\eta^{1/2} \ibeta[\big]{j_1+j_2, -j_2+\frac{1}{2}}{\eta}$.

Comparing with the result \eqref{eq:solution-l=1}, we should have $\alpha_2=0$.
This can be deduced, independently, by passing to the bulk (non-chiral) correlator and imposing single-valuedness.
This bulk correlator has the form
\begin{multline}
	\bracket[\big]{\phi_{j_1}}{\phi_{j_2}(1,1) \phi_{1/2}(\eta,\ahol{\eta})}{\phi_{j_4}^1}
	= \delta_{j_1+j_2+j_4=1/2} \abs[\big]{\eta} \biggl( \alpha_{11} + \alpha_{12}\ibeta[\big]{-j_4+\tfrac{1}{2}, -j_2+\tfrac{1}{2}}{\ahol{\eta}} \biggr. \\
	\biggl. + \alpha_{21}\ibeta[\big]{-j_4+\tfrac{1}{2}, -j_2+\tfrac{1}{2}}{\eta} + \alpha_{22}\abs[\Big]{\ibeta[\big]{-j_4+\tfrac{1}{2}, -j_2+\tfrac{1}{2}}{\eta}}^2 \biggr),
\end{multline}
where the $\alpha_{ij}$ are unknown constants.
Here, we assume that the bulk ghost primaries factorise (which holds when the charges are not integers) and that the chiral and antichiral charges agree and are real (this is the diagonal modular invariant).

To ensure single-valuedness, first consider the monodromy around $\eta=\ahol{\eta}=0$.
Under $\eta \mapsto \ee^{2\pi\ii} \eta$ and $\ahol{\eta} \mapsto \ee^{-2\pi\ii} \ahol{\eta}$, the terms involving $\alpha_{11}$ and $\alpha_{22}$ are invariant, while those involving $\alpha_{12}$ and $\alpha_{21}$ pick up factors of $-\ee^{2\pi\ii j_4}$ and $-\ee^{-2\pi\ii j_4}$, respectively, see \eqref{eq:def-Beta}.
Assuming that $j_4 \notin \ZZ+\frac{1}{2}$, we conclude that $\alpha_{12} = \alpha_{21} = 0$.
A similar analysis of the behaviour around $\eta=\ahol{\eta}=1$, using \eqref{eq:beta-0-to-1}, however concludes that $\alpha_{22} = 0$ when $j_2 \notin \ZZ+\frac{1}{2}$.
Assuming that the bulk correlator is a continuous function of the charges, the final result is thus
\begin{equation}
  \bracket[\big]{\phi_{j_1}}{\phi_{j_2}(1,1) \phi_{1/2}(\eta,\ahol{\eta})}{\phi_{j_4}^1}
  = \delta_{j_1+j_2+j_4=1/2} \alpha_{11} \abs[\big]{\eta}, \quad \text{where}\ \alpha_{11} = \abs[\Big]{\corrconst{0&0&0&1}{j_1&j_2&1/2&j_4}}^2.
\end{equation}
This is in perfect agreement with the chiral solution \eqref{eq:solution-l=1} obtained from the KZ equations \eqref{eq:kz-method1} and \eqref{eq:kz-J0-l=1}.

Similar calculations reveal the coordinate dependence of the $4$-point function \eqref{eq:special-4pt-function} with $j_3=\frac{1}{2}$ and $\ell=3$.
The analogue of the chiral correlator \eqref{eq:4pt-sol-l=1} is
\begin{equation} \label{eq:4pt-sol-l=3}
	\bracket[\big]{\phi_{j_1}}{\phi_{j_2}(1) \phi_{1/2}(\eta)}{\phi_{j_4}^3}
	= \delta_{j_1+j_2+j_4=5/2} \eta^{-j_4+2} (1-\eta)^{-j_2+1/2} \brac[\Big]{\alpha_1 + \alpha_2 \ibeta[\big]{j_4-\tfrac{1}{2}, j_2-\tfrac{1}{2}}{\eta}}
\end{equation}
and single-valuedness of the corresponding bulk correlator again rules out the incomplete beta function as a (contributing) conformal block.
We therefore have
\begin{equation}
  \bracket[\big]{\phi_{j_1}}{\phi_{j_2}(1,1) \phi_{1/2}(\eta,\ahol{\eta})}{\phi_{j_4}^3}
  = \delta_{j_1+j_2+j_4=5/2} \alpha_{11} \abs{\eta}^{-2j_4+4} \abs{1-\eta}^{-2j_2+1},
\end{equation}
which is consistent with the chiral correlator being
\begin{equation} \label{eq:4pt-l=3}
  \bracket[\big]{\phi_{j_1}}{\phi_{j_2}(1) \phi_{1/2}(\eta)}{\phi_{j_4}^3}
  = \delta_{j_1+j_2+j_4=5/2} \corrconst{0&0&0&3}{j_1&j_2&1/2&j_4} \eta^{-j_4+2} \brac{1-\eta}^{-j_2+1/2}.
\end{equation}
As in the $\ell=1$ case, there are no logarithmic singularities for $\ell=3$.

We remark that if we substitute $j_2=0$ into this solution, then we expect to get a (partially specialised) $3$-point function because $\phi_0$ should be identified with the vacuum.
However, this $3$-point function should vanish identically according to \zcref{sec:kz} and \eqref{eq:3pt-BPZ}.
This is only consistent with \eqref{eq:4pt-l=3} if the $4$-pt constant vanishes when $j_2$ does: $\corrconst{0&0&0&3}{j_1&0&1/2&j_4} = 0$.
Similarly, we must also have $\corrconst{0&0&0&3}{0&j_2&1/2&j_4} = 0$.

Our last remaining chance to find logarithmic singularities is thus the case $\ell=2$.
Now, the solutions of the BPZ equation are (fully) hypergeometric:
\begin{equation} \label{eq:4pt-sol-l=2}
	\begin{split}
		\bracket[\big]{\phi_{j_1}}{\phi_{j_2}(1) \phi_{1/2}(\eta)}{\phi_{j_4}^2}
		= \delta_{j_1+j_2+j_4=3/2} \eta \brac[\Big]{\alpha_1 \,\twoFone*{-j_1+1}{1/2}{j_4+1/2}{\eta} + \alpha_2 \eta^{-j_4+1/2} \,\twoFone*{j_2}{-j_4+1}{-j_4+3/2}{\eta}}.
	\end{split}
\end{equation}
Imposing the single-valuedness of the bulk correlator around $0$ leads this time to the following combination of conformal blocks:
\begin{equation} \label{eq:4pt-bulk-sol-l=2}
  \bracket[\big]{\phi_{j_1}}{\phi_{j_2}(1,1) \phi_{1/2}(\eta,\ahol{\eta})}{\phi_{j_4}^2}
  = \delta_{j_1+j_2+j_4=3/2} \abs{\eta}^2
  \brac[\Bigg]{\alpha_{11} \abs[\bigg]{\twoFone*{-j_1+1}{1/2}{j_4+1/2}{\eta}}^2 + \alpha_{22} \abs{\eta}^{-2j_4+1} \abs[\bigg]{\twoFone*{j_2}{-j_4+1}{-j_4+3/2}{\eta}}^2}.
\end{equation}
The single-valuedness around $1$ may be analysed by transforming both hypergeometric functions with \eqref{eq:hypergeom-0-to-1} and then applying \eqref{eq:euler-trans-formula} to the hypergeometric functions multiplying $\alpha_{22}$.
The cross-terms that arise from squaring are now seen to have non-trivial monodromy around $1$ unless their coefficients vanish.
This leads to the following constraint:
\begin{equation} \label{eq:l=2-constraint}
	\begin{split}
		\alpha_{22}
	  &= -\frac{\Gamma(j_4+\frac{1}{2})}{\Gamma(-j_1+1)\Gamma(-j_2+1)}
	      \frac{\Gamma(j_4+\frac{1}{2})}{\Gamma(\frac{1}{2})\Gamma(j_4)}
	      \frac{\Gamma(j_1)\Gamma(j_2)}{\Gamma(\frac{3}{2}-j_4)}
	      \frac{\Gamma(\frac{1}{2})\Gamma(1-j_4)}{\Gamma(\frac{3}{2}-j_4)} \alpha_{11} \\
	  &= -\frac{\cbeta{j_1,j_2} \cbeta{\tfrac{1}{2},-j_4+1}}{\cbeta{-j_1+1,-j_2+1} \cbeta{\tfrac{1}{2},j_4}} \alpha_{11}.
	\end{split}
\end{equation}
Here, $\beta$ denotes the usual beta function \eqref{eq:def-beta}.
In any case, the conclusion is that the proportionality constant is non-zero and finite if we avoid $j_1, j_2, j_4 \in \ZZ$ and $j_4 \in \ZZ+\frac{1}{2}$.

The latter condition is different to the others, which may be explained as follows.
Recall that $\twoFone*{a}{b}{c}{\eta}$ and $\eta^{1-c} \,\twoFone*{a-c+1}{b-c+1}{2-c}{\eta}$ are independent solutions of the hypergeometric differential equation if and only if $c \notin \ZZ$.
The conformal blocks of \eqref{eq:4pt-sol-l=2} are therefore independent if and only if $j_4 \notin \ZZ+\frac{1}{2}$.
When $j_4 \in \ZZ+\frac{1}{2}$, one of these hypergeometric functions must be replaced by a more complicated solution, namely $\log \eta$ times the other hypergeometric function plus some other series in $\eta$.
A more involved study of single-valuedness then shows that the ratio of the analogues of $\alpha_{11}$ and $\alpha_{22}$ is also non-zero and finite when $j_4 \in \ZZ+\frac{1}{2}$.

Incidentally, the case when at least one of the charges is an integer corresponds to the corresponding module being reducible.
It is therefore no surprise that the correlation function analysis is much more delicate for these cases.

Our conclusion is then that the bulk $\ell=2$ correlator \eqref{eq:4pt-bulk-sol-l=2} exhibits logarithmic singularities at $\eta=\ahol{\eta}=0$ when $j_4 \in \ZZ+\frac{1}{2}$.
Since $j_3$ is fixed to $\frac{1}{2}$, these singularities arise when $j_3+j_4 \in \ZZ$, which in turn requires that $j_1+j_2 \in \ZZ$ (for the correlator to be non-vanishing).
From the perspective of the fusion rules \eqref{eq:fusion-rule-2}, $j_1+j_2 \in \ZZ$ implies that $\relmod{j_1}$ and $\relmod{j_2}$ fuse to give $\sfproj{-1}$, while $j_3+j_4 \in \ZZ$ implies that $\relmod{j_3}$ and $\sfrel{2}{j_4}$ fuse to give $\sfproj{}$.
This appearance of logarithmic singularities is thus in perfect agreement with the intuition that they arise because these \opes\ generate two Jordan partner fields that are conjugate to one another.

Applying the transformation formula \eqref{eq:hypergeom-0-to-1} to our conformal blocks, a comparison with \eqref{eq:hypergeom-around-1} shows that the bulk $\ell=2$ correlator will also develop logarithmic singularities at $\eta = \ahol{\eta} = 1$ when $a+b-c+1 \in \ZZ$, hence $j_2 \in \ZZ+\frac{1}{2}$.
In this case, we have $j_2+j_3 \in \ZZ$, hence $j_1+j_4 \in \ZZ$, and the interpretation in terms of fusion rules is that now $\relmod{j_2}$ and $\relmod{j_3}$ fuse to give $\sfproj{-1}$, while $\relmod{j_1}$ and $\sfrel{2}{j_4}$ fuse to give $\sfproj{}$.
Again, two conjugate Jordan partner fields are generated, explaining the logarithmic singularity.

We conclude by noting that if we set $j_2$ to $0$ in the bulk correlator \eqref{eq:4pt-bulk-sol-l=2}, then the conformal blocks become
\begin{equation}
  \eta \, \twoFone*{j_4-1/2}{1/2}{j_4+1/2}{\eta} \quad \text{and} \quad \eta^{-j_4+3/2} \, \twoFone*{0}{-j_4+1}{-j_4+1/2}{\eta} = \eta^{j_1},
\end{equation}
while the appearance of $\Gamma(0)$ in the numerator of \eqref{eq:l=2-constraint} indicates that we must have $\alpha_{11} = 0$.
The result is thus
\begin{equation}
  \bracket[\big]{\phi_{j_1}}{\phi_{1/2}(\eta,\ahol{\eta})}{\phi_{j_4}^2} = \delta_{j_1+j_4=3/2} \alpha_{22} \abs{\eta}^{2j_1},
\end{equation}
which agrees with the (partially specialised) $3$-point function obtained from \eqref{eq:special-3pt-function} if we set $\alpha_{22} = \abs[\Big]{\corrconst{0&0&2}{j_1&1/2&j_4}}^2$.
It is easy to check that \eqref{eq:4pt-bulk-sol-l=2} likewise agrees with \eqref{eq:special-3pt-function} when we set $j_1$ to $0$.

\section{Further explorations} \label{sec:exploring}

We have computed the ghost-primary $4$-point functions \eqref{eq:special-4pt-function} when $\ell=1$ (for arbitrary charges) and when $\ell=2$ or $3$ (with the restriction that $j_3=\frac{1}{2}$).
It is therefore natural to ask, in the latter case, if we can obtain the $4$-point functions when $j_3$ takes other values.
With the tools developed here, a general solution appears to be out of reach.
However, we can answer the question for $j_3 \in \NN+\frac{1}{2}$ by inputting a (specialised) $j_3 = \frac{1}{2}$ conformal block into the following algorithm.
\begin{itemize}
  \item First, determine the function $H(\eta)$ from the following version of \eqref{eq:limiting-N-point-functions}:
  \begin{equation}
    \bracket[\big]{\phi_{j_1}}{\phi_{j_2}(1) \phi_{j_3}(\eta)}{\phi_{j_4}^{\ell}}
		= \delta_{j_1+j_2+j_3+j_4=\ell} \eta^{-2h_4^{\ell} + (j_3+j_4^{\ell})/2} (1-\eta)^{h_4^{\ell}/3+(j_2+j_3)/2} H(\eta).
  \end{equation}
  Here, we set $j_4^{\ell} = j_4-\ell$ and $h_4^{\ell} = j_4 \ell - \frac{1}{2} \ell(\ell+1)$, for convenience.
  \item Now, substitute this into the following version of \eqref{eq:Ward-identity-sols}:
  \begin{multline}
    \corrfn[\big]{\phi_{j_1}(w_1) \phi_{j_2}(w_2) \phi_{j_3}(w_3) \phi_{j_4}^{\ell}(w_4)}
    = \delta_{j_1+j_2+j_3+j_4=\ell} H(\eta) \\
      \cdot w_{12}^{h_4^{\ell}/3 + (j_1+j_2)/2} w_{13}^{h_4^{\ell}/3 + (j_1+j_3)/2} w_{23}^{h_4^{\ell}/3 + (j_2+j_3)/2}
      w_{14}^{-2h_4^{\ell}/3 + (j_1+j_4^{\ell})/2} w_{24}^{-2h_4^{\ell}/3 + (j_2+j_4^{\ell})/2} w_{34}^{-2h_4^{\ell}/3 + (j_3+j_4^{\ell})/2}.
  \end{multline}
  This is the unspecialised conformal block.
  \item Next, we substitute into the KZ equation \eqref{eq:kz-method2} with $N=4$, rearranged into the form
  \begin{multline} \label{eq:kz-method2-bpz}
    \corrfn[\big]{\phi_{j_1}(w_1) \phi_{j_2}(w_2) \phi_{j_3+1}(w_3) \phi_{j_4-1}^{\ell}(w_N)} \\
    = -\frac{w_{34}^{\ell+1}}{j_3} \brac[\bigg]{\pd_4 + (\ell-1) \brac[\Big]{\frac{j_1}{w_{14}} + \frac{j_2}{w_{24}} + \frac{j_3}{w_{34}}}} \corrfn[\big]{\phi_{j_1}(w_1) \phi_{j_2}(w_2) \phi_{j_3}(w_3) \phi_{j_4}^{\ell}(w_N)} \\
    - \frac{j_1}{j_3} \brac*{\frac{w_{34}}{w_{14}}}^{\ell+1} \corrfn[\big]{\phi_{j_1+1}(w_1) \phi_{j_2}(w_2) \phi_{j_3}(w_3) \phi_{j_4-1}^{\ell}(w_N)}
    - \frac{j_2}{j_3} \brac*{\frac{w_{34}}{w_{24}}}^{\ell+1} \corrfn[\big]{\phi_{j_1}(w_1) \phi_{j_2+1}(w_2) \phi_{j_3}(w_3) \phi_{j_4-1}^{\ell}(w_N)}.
  \end{multline}
  \item Finally, specialise to $w_1 \to \infty$, $w_2 = 1$, $w_3 = \eta$ and $w_4 = 0$.
  The result is
  \begin{equation*}
    \bracket[\big]{\phi_{j_1}}{\phi_{j_2}(1) \phi_{j_3+1}(\eta)}{\phi_{j_4-1}^{\ell}}.
  \end{equation*}
\end{itemize}

We test this algorithm by verifying that it reproduces our $\ell=1$ results.
Indeed, applying it to the specialised $j_3=\frac{1}{2}$ correlator
\begin{equation}
  \bracket{\phi_{j_1}}{\phi_{j_2}(1) \phi_{1/2}(w_3)}{\phi_{j_4}^1}
  = \delta_{j_1+j_2+j_4=1/2} \corrconst{0&0&0&1}{j_1&j_2&1/2&j_4} \eta^{1/2}
\end{equation}
gives the corresponding $j_3=\frac{3}{2}$ correlator:
\begin{equation}
  \bracket{\phi_{j_1}}{\phi_{j_2}(1) \phi_{3/2}(w_3)}{\phi_{j_4-1}^1}
  = \delta_{j_1+j_2+j_4=1/2} \corrconst{0&0&0&1}{j_1&j_2&1/2&j_4} \eta^{3/2}.
\end{equation}
Substituting the result into the algorithm iteratively, we obtain
\begin{equation}
  \bracket{\phi_{j_1}}{\phi_{j_2}(1) \phi_{1/2+k}(w_3)}{\phi_{j_4-k}^1}
  = \delta_{j_1+j_2+j_4=1/2} \corrconst{0&0&0&1}{j_1&j_2&1/2&j_4} \eta^{1/2+k}, \quad k \in \NN.
\end{equation}
This is in full agreement with the general result \eqref{eq:solution-l=1} provided that we identify $\corrconst{0&0&0&1}{j_1&j_2&k+1/2&j_4-k}$ with $\corrconst{0&0&0&1}{j_1&j_2&1/2&j_4}$, for all $k \in \NN$.
(The latter identification is of course also implied by \eqref{eq:4pt-constants-l=1}.)
Moreover, it suggests that we can even interpolate from $j_3 \in \NN+\frac{1}{2}$ to correctly guess the general $j_3 \in \RR$ result.

With this successful test, we next apply the algorithm to the $\ell=3$ specialised correlator \eqref{eq:4pt-l=3}, for which the results will be new.
We obtain (again for $k \in \NN$)
\begin{equation} \label{eq:4pt-block-l=3}
  \bracket{\phi_{j_1}}{\phi_{j_2}(1) \phi_{1/2+k}(\eta)}{\phi_{j_4-k}^3}
  = \delta_{j_1+j_2+j_4=5/2} \corrconst{0&0&0&3}{j_1&j_2&1/2&j_4} \eta^{-j_4+3k+2} (1-\eta)^{-j_2+1/2} P_k(j_1,j_4;\eta),
\end{equation}
where $P_k$ is the following polynomial in $\eta$:
\begin{equation}
  P_k(j_1,j_4;\eta) = (-1)^k \frac{2^k}{(2k-1)!!} \sum_{i=0}^k \binom{k}{i} (-j_1+\tfrac{3}{2})_i (j_4-k)_{k-i} \eta^i.
\end{equation}
Here, $(a)_n$ denotes the Pochhammer symbol, see \eqref{eq:def-hypergeom}.
This does not obviously generalise to $k \notin \NN$ until we note that
\begin{equation} \label{eq:poly-to-hyper-l=3}
  P_k(j_1,j_4;\eta)
  = \frac{\Gamma(\frac{1}{2})}{\Gamma(k+\frac{1}{2})} \frac{\Gamma(-j_4+k+1)}{\Gamma(-j_4+1)} \,\twoFone*{-j_1+3/2}{-k}{-j_4+1}{\eta}
  = \frac{\cbeta{\tfrac{1}{2},-j_4+k+1}}{\cbeta{k+\tfrac{1}{2},-j_4+1}} \,\twoFone*{-j_1+3/2}{-k}{-j_4+1}{\eta}.
\end{equation}
It is therefore not unreasonable to conjecture that \eqref{eq:4pt-block-l=3} is correct for all $k \in \RR$, hence all $j_3 \in \RR$, once we replace the polynomial $P_k$ by its hypergeometric counterpart \eqref{eq:poly-to-hyper-l=3}:
\begin{multline} \label{eq:conj-4pt-block-l=3}
	\bracket{\phi_{j_1}}{\phi_{j_2}(1) \phi_{j_3}(\eta)}{\phi_{j_4}^3} = \delta_{j_1+j_2+j_3+j_4=3} \corrconst{0&0&0&3}{j_1&j_2&1/2&j_3+j_4-1/2} \\
	\cdot \frac{\cbeta{\frac{1}{2},-j_4+1}}{\cbeta{j_3,-j_3-j_4+\frac{3}{2}}} \eta^{2j_3-j_4+1} (1-\eta)^{-j_2+1/2} \twoFone*{-j_1+3/2}{-j_3+1/2}{-j_3-j_4+3/2}{\eta}.
\end{multline}
This also suggests that the general $4$-point constant $\corrconst{0&0&0&3}{j_1&j_2&j_3&j_4}$ is the product of the $\Gamma$-ratios and the $4$-point constant in \eqref{eq:conj-4pt-block-l=3}.

As a sanity check, consider the $4$-point function \eqref{eq:4pt-block-l=3} with $j_3=-\frac{1}{2}$.
According to \eqref{eq:conj-4pt-block-l=3}, it takes the form
\begin{equation} \label{eq:4-pt-block-l=3-j3=-1/2}
	\bracket{\phi_{j_1}}{\phi_{j_2}(1) \phi_{-1/2}(\eta)}{\phi_{j_4}^3} = \delta_{j_1+j_2+j_4=7/2} \frac{\corrconst{0&0&0&3}{j_1&j_2&1/2&j_4-1}}{2(j_4-1)}
	\eta^{-j_4} (1-\eta)^{-j_2+1/2} \,\twoFone*{-j_1+3/2}{1}{-j_4+2}{\eta}.
\end{equation}
Substituting this into our algorithm then indeed returns the correct $j_3=+\frac{1}{2}$ correlator \eqref{eq:4pt-l=3}.
One can similarly verify that similar checks pass for $k=-2,-3,\dots$, showing that \eqref{eq:conj-4pt-block-l=3} is correct for $j_3 \in \ZZ+\frac{1}{2}$.

Note that the factor $\cbeta{j_3,-j_3-j_4+\frac{3}{2}}$ in the denominator of \eqref{eq:conj-4pt-block-l=3} involves a factor $\Gamma(-j_3-j_4+3/2)$ that regularises the hypergeometric function in the numerator.
The ratio is thus defined even when $-j_3-j_4+3/2 \in -\NN$.
In particular, this shows that \eqref{eq:conj-4pt-block-l=3} only admits power-law, and not logarithmic, singularities at $\eta=0$.
However, the hypergeometric function will generically admit a logarithmic singularity at $\eta=1$, if $j_1-j_4 \in \ZZ+\frac{1}{2}$, and at $\eta=\infty$, if $j_1-j_3 \in \ZZ$.
Neither of these singularities are regularised in \eqref{eq:conj-4pt-block-l=3} and these charge conditions are inconsistent with those producing logarithmic modules through fusion.
We therefore conclude that the validity of \eqref{eq:conj-4pt-block-l=3} requires that $\corrconst{0&0&0&3}{j_1&j_2&1/2&j_3+j_4-1/2}$ vanishes when $j_1-j_4 \in \ZZ+\frac{1}{2}$ or $j_1-j_3 \in \ZZ$.

Another curious feature of \eqref{eq:conj-4pt-block-l=3} is that there may exist divergences as the charges take special values.
This is already evident in the $j_3=-\frac{1}{2}$ correlator \eqref{eq:4-pt-block-l=3-j3=-1/2}.
Taking $j_4 \to 1$ sends the hypergeometric function to $(1-\eta)^{j_1-3/2}$, so the only way to avoid a divergence is to have $\corrconst{0&0&0&3}{j_1&j_2&1/2&0} = 0$.
This is not unrealistic given that we deduced in \zcref{sec:bpz} analogous vanishing results when $j_1$ or $j_2$ is $0$.

Our algorithm returns even more tantalising results for the $\ell=2$ specialised correlators.
We start with the conformal blocks of \eqref{eq:4pt-sol-l=2}:
\begin{equation}
	\begin{aligned}
		\confblock{1}{0&0&0&2}{j_1&j_2&1/2&j_4}{\eta} &= \delta_{j_1+j_2+j_4=3/2} \eta \,\twoFone*{-j_1+1}{1/2}{j_4+1/2}{\eta} \\
		\text{and} \quad
		\confblock{2}{0&0&0&2}{j_1&j_2&1/2&j_4}{\eta} &= \delta_{j_1+j_2+j_4=3/2} \eta^{-j_4+3/2} \,\twoFone*{j_2}{-j_4+1}{-j_4+3/2}{\eta}.
	\end{aligned}
\end{equation}
Using these as inputs in the algorithm then results in conformal blocks for which $j_3 = \frac{1}{2}$ is replaced by $k+\frac{1}{2}$ and $j_4$ by $j_4-k$, for $k \in \NN$.
In this way, we obtain
\begin{equation}
	\begin{aligned}
		\confblock{1}{0&0&0&2}{j_1&j_2&k+1/2&j_4-k}{\eta} &= \delta_{j_1+j_2+j_4=3/2} \eta^{2k+1}
		\sum_{m=0}^k \binom{k}{m} \frac{(-j_4+\frac{1}{2})_{k-m} (j_4)_m}{(\frac{1}{2})_k} \,\twoFone*{-m+1/2}{-j_1+1}{j_4+1/2}{\eta} \\
		\text{and} \quad
		\confblock{2}{0&0&0&2}{j_1&j_2&k+1/2&j_4-k}{\eta} &= \delta_{j_1+j_2+j_4=3/2} \eta^{2k-j_4+3/2}
		\sum_{m=0}^k \binom{k}{m} \frac{(-j_4+\frac{1}{2})_{k-m} (\frac{1}{2})_m}{(\frac{1}{2})_k} \,\twoFone*{-j_4-m+1}{j_2}{-j_4+3/2}{\eta}.
	\end{aligned}
\end{equation}
Here, we have used \eqref{eq:diff-hypergeom} and \eqref{eq:contiguous} liberally in simplifying the hypergeometric functions.
We note that monodromy considerations show that the coefficients of the squares of these blocks are still given by \eqref{eq:l=2-constraint} in the bulk correlator, independent of $k \in \NN$.

Remarkably, these sums may be evaluated in closed form so as to suggest generalisations to all $k \in \RR$:
\begin{equation}
	\begin{aligned} \label{eq:l2-conformal-blocks-sum-form}
		\sum_{m=0}^k \binom{k}{m} \frac{(-j_4+\frac{1}{2})_{k-m} (j_4)_m}{(\frac{1}{2})_k} \,\twoFone*{-m+1/2}{-j_1+1}{j_4+1/2}{\eta}
		&= (1-\eta)^{j_1-1} \,\threeFtwo[\Big]{j_4}{-j_1+1}{k+1/2}{j_4+1/2}{1/2}{-\frac{\eta}{1-\eta}} \\
		\text{and} \quad
		\sum_{m=0}^k \binom{k}{m} \frac{(-j_4+\frac{1}{2})_{k-m} (\frac{1}{2})_m}{(\frac{1}{2})_k} \,\twoFone*{-j_4-m+1}{j_2}{-j_4+3/2}{\eta}
		&= \frac{(-j_4+1)_k}{(\frac{1}{2})_k} (1-\eta)^{-j_2} \,\threeFtwo[\Big]{1/2}{j_2}{-j_4+k+1}{-j_4+3/2}{-j_4+1}{-\frac{\eta}{1-\eta}}.
	\end{aligned}
\end{equation}
We recall the (series) definition of the generalised hypergeometric function in \eqref{eq:def-3F2}.
A derivation of these summation formulae is also provided in \zcref{sec:derivation}, for convenience.
The corresponding conjectures for the conformal blocks are thus
\begin{equation} \label{eq:conj-4pt-blocks-l=2}
	\begin{aligned}
		\confblock{1}{0&0&0&2}{j_1&j_2&j_3&j_4}{\eta} &= \delta_{j_1+j_2+j_3+j_4=2} \eta^{2j_3}
		(1-\eta)^{j_1-1} \,\threeFtwo[\Big]{j_3+j_4-1/2}{-j_1+1}{j_3}{j_3+j_4}{1/2}{-\frac{\eta}{1-\eta}} \\
		\text{and} \quad
		\confblock{2}{0&0&0&2}{j_1&j_2&j_3&j_4}{\eta} &= \delta_{j_1+j_2+j_3+j_4=2} \frac{\cbeta{\frac{1}{2}, j_3-j_4+\frac{1}{2}}}{\cbeta{j_3,1-j_4}}
		\eta^{j_3-j_4+1} (1-\eta)^{-j_2} \,\threeFtwo[\Big]{1/2}{j_2}{-j_4+1}{j_1+j_2}{j_1+j_2-1/2}{-\frac{\eta}{1-\eta}}.
	\end{aligned}
\end{equation}

Finally, we offer some speculation as to why the ``interpolation conjectures'' \eqref{eq:conj-4pt-block-l=3} and \eqref{eq:conj-4pt-blocks-l=2} are not completely unreasonable, even though they generalise from $j_3 \in \NN+\frac{1}{2}$ to $j_3 \in \RR$.
For this, observe that the ghost modules $\relmod{j}$, for $0<j<1$, and $\relmod{0}^+$ are almost identical.
Aside from the reducibility of $\relmod{0}^+$, the only real difference is that the charges lie in $\ZZ+j$.
This can be formalised \cite{MatCla00,MazLec10} by saying that these modules form a \emph{coherent family}.

Rather than dwell on a precise definition, we instead note that every member of this family can be obtained from any other member by applying functors known as \emph{twisted localisation}.
These functors are defined as follows:
\begin{itemize}
	\item First, we localise the mode algebra of the bosonic ghosts \voa\ $\bgvoa$ with respect to $\gamma_0$.
	The result is the extension of the mode algebra by the formal inverse $\gamma_0^{-1}$.
	We shall call it the \emph{localised algebra}.
	It satisfies
	\begin{equation}
		\gamma_0^{-1} \gamma_0 = \wun = \gamma_0 \gamma_0^{-1}, \quad
		\comm{\beta_m}{\gamma_0^{-1}} = \delta_{m=0} \gamma_0^{-2} \quad \text{and} \quad
		\comm{\gamma_m}{\gamma_0^{-1}} = 0,
	\end{equation}
	where $\gamma_0^{-2}$ is obviously denoting $(\gamma_0^{-1})^2$.
	These relations are motivated by computing $\comm{\beta_m}{\gamma_0^k} = -k \gamma_0^{k-1}$, for all $k \in \NN$, \emph{and then setting $k=-1$}.
	It is easy to check that these relations make the localised algebra associative.
	\item Next, we induce some $\bgvoa$-module to obtain a module for the localised algebra.
	In general, this annihilates states when $\gamma_0$ does not act injectively and creates states when $\gamma_0$ acts injectively but not bijectively.
	However, $\gamma_0$ acts bijectively on all the modules in our chosen coherent family.
	Inducing these modules therefore does nothing except allow us to act with $\gamma_0^{-1}$.
	\item The localised algebra admits automorphisms $\Theta_k$, $k \in \ZZ$, defined by conjugating with $\gamma_0^k$.
	It is easy to check that $\Theta_k$ preserves $\gamma_0^{-1}$ and all the ghost modes \emph{except} $\beta_0$.
	For the latter, we find that
	\begin{equation}
		\Theta_k(\beta_0) = \gamma_0^k \beta_0 \gamma_0^{-k} = \gamma_0^k \gamma_0^{-k} \beta_0 + \gamma_0^k \comm{\beta_0}{\gamma_0^{-k}}
		= \beta_0 + k \gamma_0^{-1}.
	\end{equation}
	But, this action now makes sense \emph{for arbitrary $k \in \RR$} (in fact, for $k \in \CC$) and it is easy to check that these generalised $\Theta_k$, $k \in \RR$, are still automorphisms of the localised algebra.
	\item It follows that we can \emph{twist} the action of the localised algebra on the induced modules by replacing each mode with its image under $\Theta_k$, $k \in \RR$.
	This will produce a new ``twisted'' module for the localised algebra.
	The key is now what this twisting does to charges and conformal weights.
	Since
	\begin{equation}
		\Theta_k(J_0) = J_0 + k \wun \quad \text{and} \quad \Theta_k(L_0) = L_0,
	\end{equation}
	all charges are shifted uniformly by $k \in \RR$ while the conformal weights are preserved.
	\item Finally, we can forget all about the action of $\gamma_0^{-1}$ by restricting our twisted module so that only the ghost modes are allowed to act.
	The result is a $\bgvoa$-module that looks identical to the one we started with, except that its charges have been uniformly shifted by $k$.
	It is not hard to verify that this twisted localisation procedure sends $\relmod{1/2}$ (say) to $\relmod{k+1/2}$ (with a superscript ``$+$'' if $k \in \ZZ+\frac{1}{2}$).
\end{itemize}

To summarise, $\relmod{0}^+$ and the $\relmod{j}$, $0<j<1$, can all be constructed from $\relmod{1/2}$ using twisted localisation functors.
The same is obviously true if we include spectral flow (now one instead localises with respect to $\gamma_{\ell}$).
Because these functors are defined by deriving relations that hold for $k \in \NN$ (or $\ZZ$) and \emph{interpolating} them to all $k \in \RR$, it seems natural to expect that all the data associated to these modules, for example correlation functions, may also be computed by interpolating the $k \in \NN$ results to all $k \in \RR$.
We shall leave this speculation here for now, noting only that we hope in the future to have data that confirms (or refutes!) this expectation.

\appendix

\section{Hypergeometric functions} \label{sec:hypergeom}

In this \lcnamezcref{sec:hypergeom}, we collect a few useful (and very well-known) facts about special functions, in particular the incomplete beta function and the hypergeometric function.

First recall that the hypergeometric differential equation is given by
\begin{equation} \label{eq:hypergeom-eqn}
	z(1-z) \frac{\dd^2y}{\dd z^2} + \brac[\big]{c-(a+b+1)z} \frac{\dd y}{\dd z} - ab y = 0,
\end{equation}
where $a$, $b$, and $c$ are complex parameters.
This equation has three regular singular points at $z=0$, $z=1$, and $z=\infty$.

Around $z=0$, there are generically two series solutions that are traditionally expressed in terms of the hypergeometric function
\begin{equation} \label{eq:def-hypergeom}
	\twoFone*{a}{b}{c}{z} = \sum_{n\ge0} \frac{(a)_n (b)_n}{(c)_n} \frac{z^n}{n!},
\end{equation}
where $(x)_n = x(x+1)\cdots(x+n-1)$ is the Pochhammer symbol (rising factorial).
This series is defined for $c \notin -\NN$ and is convergent for $\abs{z}<1$.
It admits an analytic continuation to the entire complex plane, though $1$ and $\infty$ are branch points whenever $c-a-b \notin \ZZ$.

The following identity is a straightforward consequence of the definition \eqref{eq:def-hypergeom}:
\begin{equation} \label{eq:diff-hypergeom}
	\frac{\dd}{\dd z} \,\twoFone*{a}{b}{c}{z} = \frac{ab}{c} \,\twoFone*{a+1}{b+1}{c+1}{z}.
\end{equation}
It easily follows that \eqref{eq:def-hypergeom} is a solution of \eqref{eq:hypergeom-eqn}.
Indeed, the two linearly independent solutions of \eqref{eq:hypergeom-eqn} are
\begin{equation} \label{eq:hypergeom-around-0}
	y_1(z) = \twoFone*{a}{b}{c}{z} \quad \text{and} \quad
	y_2(z) = z^{1-c} \,\twoFone*{a-c+1}{b-c+1}{2-c}{z},
\end{equation}
provided that $c \notin \ZZ$.
In the case that $c \in \ZZ$, only one series solution can be found.
The other linearly independent solution is then obtained by multiplying this solution by $\log z$ and adding another series (this is the method of Frobenius).

There are also generically independent series solutions around $z=1$ and $z=\infty$.
In particular, around $z=1$, they are
\begin{equation} \label{eq:hypergeom-around-1}
	\twoFone*{a}{b}{a+b-c+1}{1-z} \quad \text{and} \quad (1-z)^{c-a-b} \,\twoFone*{c-a}{c-b}{c-a-b+1}{1-z},
\end{equation}
assuming that $a+b-c \notin \ZZ$, while around $z=\infty$ they have the forms
\begin{equation} \label{eq:hypergeom-around-inf}
	z^{-a} \,\twoFone*{a}{a-c+1}{a-b+1}{z^{-1}} \quad \text{and} \quad z^{-b} \,\twoFone*{b}{b-c+1}{b-a+1}{z^{-1}},
\end{equation}
assuming that $a-b \notin \ZZ$.
Again, if these conditions on $a$, $b$ and $c$ are not met, then one solution will develop a logarithmic singularity.

Aside from the obvious symmetry $\twoFone*{a}{b}{c}{z} = \twoFone*{b}{a}{c}{z}$, the hypergeometric function also satisfies the following identities:
\begin{subequations}
	\begin{align}
		\twoFone*{a}{b}{c}{z} &= (1-z)^{-b} \twoFone*{c-a}{b}{c}{-\frac{z}{1-z}} & &\text{(Pfaff's transformation)}, \label{eq:Pfaff-trans-formula} \\
		\twoFone*{a}{b}{c}{z} &= (1-z)^{c-a-b} \,\twoFone*{c-a}{c-b}{c}{z} & &\text{(Euler's transformation)}. \label{eq:euler-trans-formula}
	\end{align}
\end{subequations}
Another more complicated transformation formula, relating the expansions around $z=0$ and $z=1$, will also turn out to be useful:
\begin{equation} \label{eq:hypergeom-0-to-1}
	\twoFone*{a}{b}{c}{z} = \frac{\Gamma(c)\Gamma(c-a-b)}{\Gamma(c-a)\Gamma(c-b)} \,\twoFone*{a}{b}{a+b-c+1}{1-z}
	+ \frac{\Gamma(c)\Gamma(a+b-c)}{\Gamma(a)\Gamma(b)} (1-z)^{c-a-b} \,\twoFone*{c-a}{c-b}{c-a-b+1}{1-z}.
\end{equation}
Finally, we will also make use of Gauss' contiguous relations:
\begin{equation} \label{eq:contiguous}
	\begin{aligned}
		a \brac[\Big]{\twoFone*{a+1}{b}{c}{z} - \twoFone*{a}{b}{c}{z}}
		&= \frac{(c-a) \,\twoFone*{a-1}{b}{c}{z} + (a+bz-c) \,\twoFone*{a}{b}{c}{z}}{1-z} \\ =
		b \brac[\Big]{\twoFone*{a}{b+1}{c}{z} - \twoFone*{a}{b}{c}{z}}
		&= \frac{(c-b) \,\twoFone*{a}{b-1}{c}{z} + (az+b-c) \,\twoFone*{a}{b}{c}{z}}{1-z} \\ =
		(c-1) \brac[\Big]{\twoFone*{a}{b}{c-1}{z} - \twoFone*{a}{b}{c}{z}}
		&= \frac{z}{c} \frac{(c-a)(c-b) \,\twoFone*{a}{b}{c+1}{z} + c(a+b-c) \,\twoFone*{a}{b}{c}{z}}{1-z}.
	\end{aligned}
\end{equation}

The incomplete beta function may now be defined as follows:
\begin{equation} \label{eq:def-Beta}
	\ibeta{a,b}{z} = \frac{z^a}{a} \,\twoFone*{a}{1-b}{a+1}{z}.
\end{equation}
This requires that $a \notin \ZZ_{\le0}$.
Moreover, setting $z=1$ recovers the usual (complete) beta function
\begin{equation} \label{eq:def-beta}
	\cbeta{a,b} = \twoFone*{a}{1-b}{a+1}{1} = \frac{\Gamma(a)\Gamma(b)}{\Gamma(a+b)}.
\end{equation}
For our purposes, we quote the transformation formula obtained by substituting \eqref{eq:hypergeom-0-to-1} into this definition:
\begin{equation} \label{eq:beta-0-to-1}
	\ibeta{a,b}{z} = \cbeta{a,b} - \frac{z^a (1-z)^b}{b} \,\twoFone*{1}{a+b}{b+1}{1-z}.
\end{equation}
Here, we have used the identity $\twoFone*{a}{b}{b}{z} = (1-z)^{-a}$.

\section{A generalised hypergeometric identity} \label{sec:derivation}

We derive a closed form for sums of the form
\begin{equation} \label{eq:the-sum}
	\sum_{m=0}^{k} \binom{k}{m} (\alpha)_{k-m} (c-a)_m \,\twoFone*{a-m}{b}{c}{\eta}.
\end{equation}
Such sums appear in the $\ell=2$ conformal blocks \eqref{eq:l2-conformal-blocks-sum-form} with $\alpha = 1-c$ or $c-1$.
To accomplish this derivation, we manipulate \eqref{eq:the-sum} as follows:
\begin{itemize}
	\item First, apply \zcref{eq:Pfaff-trans-formula} (Pfaff's transformation) to obtain
	\begin{equation} \label{eqn:apdx-step1-pfaff-transform}
		(1-\eta)^{-b} \sum_{m=0}^{k} \binom{k}{m} (\alpha)_{k-m} (c-a)_m \,\twoFone*{c-a+m}{b}{c}{w}, \quad
		\text{where}\ w = -\frac{\eta}{1-\eta}.
	\end{equation}
	\item Next, we rearrange \cite[Eq.~15.5.3]{NIST_DLMF_2010} to obtain
	\begin{equation} \label{eqn:apdx-step2-mth-differential}
		(c-a)_m \,\twoFone*{c-a+m}{b}{c}{w} = w^{a-c-m+1} \brac*{w \frac{\dd}{\dd w} w}^m \brac*{w^{c-a-1} \twoFone*{c-a}{b}{c}{w}}.
	\end{equation}
	\item Letting $z=w^{-1}$, we note that
	\begin{equation}
		w \frac{\dd}{\dd w} w = -z \frac{\dd}{\dd z} z^{-1} \quad \Ra \quad
		\brac*{w \frac{\dd}{\dd w} w}^m = (-1)^m \brac*{z \frac{\dd^m}{\dd z^m} z^{-1}}.
	\end{equation}
	\eqref{eqn:apdx-step1-pfaff-transform} thus takes the form
	\begin{equation} \label{eqn:apdx-step2-2nd-3rd-factor}
		(1-\eta)^{-b} \sum_{m=0}^{k} \binom{k}{m} (\alpha)_{k-m} (-1)^m z^{c-a+m} \frac{\dd^m}{\dd z^m}\left( z^{a-c} \twoFone*{c-a}{b}{c}{z^{-1}} \right).
	\end{equation}
	\item Now write $(\alpha)_{k-m}$ as a $(k-m)$-th derivative:
	\begin{equation} \label{eqn:apdx-step3-Pochhammer-in-terms-derivative}
		(\alpha)_{k-m} = (-1)^{k-m} z^{\alpha+k-m} \frac{\dd^{k-m}}{\dd z^{k-m}} z^{-\alpha}.
	\end{equation}
	Substituting this into \eqref{eqn:apdx-step2-2nd-3rd-factor} allows us to use the Leibniz rule.
	The result is
	\begin{equation}
		\begin{split}
			&(-1)^k (1-\eta)^{-b} z^{\alpha-a+c+k} \sum_{m=0}^{k} \binom{k}{m} \frac{\dd^{k-m}}{\dd z^{k-m}} \brac[\big]{z^{-\alpha}} \frac{\dd^m}{\dd z^m} \brac*{z^{a-c} \twoFone*{c-a}{b}{c}{z^{-1}}} \\
			={} &(-1)^k (1-\eta)^{-b} z^{\beta+k} \frac{\dd^k}{\dd z^k} \brac*{z^{-\beta} \twoFone*{c-a}{b}{c}{z^{-1}}},
		\end{split}
	\end{equation}
	where $\beta=\alpha-a+c$.
	\item Insert the series definition \eqref{eq:def-hypergeom} of the hypergeometric function and perform the differentiation.
	We arrive at
	\begin{equation}
		\begin{split}
			&(-1)^k (1-\eta)^{-b} \sum_{n\ge0} \frac{(c-a)_n (b)_n}{(c)_n} (-\beta-n)(-\beta-n-1)\cdots(-\beta-n-k+1) \frac{z^{-n}}{n!} \\
			={} &(1-\eta)^{-b} \sum_{n\ge0} \frac{(c-a)_n (b)_n}{(c)_n} (\beta+n)_k \frac{z^{-n}}{n!}.
		\end{split}
	\end{equation}
	\item It only remains to massage $(\beta+n)_k$ using the identity
	\begin{equation}
		(\beta)_k (\beta+k)_n = (\beta)_n (\beta+n)_k.
	\end{equation}
	The resulting sum is then readily recognised as a generalised hypergeometric function of the form
	\begin{equation} \label{eq:def-3F2}
		\threeFtwo*{a}{b}{c}{d}{e}{z} = \sum_{n\ge0} \frac{(a)_n (b)_n (c)_n}{(d)_n (e)_n} \frac{z^n}{n!}.
	\end{equation}
\end{itemize}
Substituting for $\beta$ and $z$, this proves the desired identity:
\begin{equation} \label{eq:the-identity}
	\sum_{m=0}^{k} \binom{k}{m} (\alpha)_{k-m} (c-a)_m \,\twoFone*{a-m}{b}{c}{\eta}
	= (\alpha-a+c)_k (1-\eta)^{-b} \threeFtwo[\Big]{c-a}{b}{\alpha-a+c+k}{c}{\alpha-a+c}{-\frac{\eta}{1-\eta}}.
\end{equation}
We remark that one can also derive this identity by starting with either \cite[Eq.~16.10.2]{NIST_DLMF_2010} or \cite[Eq.~8]{karlsson_hypergeometric_1971}.

\flushleft
\providecommand{\opp}[2]{\textsf{arXiv:\mbox{#2}/#1}}
\providecommand{\pp}[2]{\textsf{arXiv:#1 [\mbox{#2}]}}

\end{document}